\documentclass[12pt,reqno]{amsart}
\usepackage[headings]{fullpage}
\usepackage{amssymb,amsmath,mathtools,bbm,tikz,tikz-cd,color,relsize,multirow}
\usepackage{bm}
\usepackage[all,cmtip]{xy}
\usepackage{url}
\usepackage{soul}
\usepackage[group-separator={,},group-minimum-digits={3}]{siunitx}
\usetikzlibrary{positioning}
\usetikzlibrary{calc}
\usetikzlibrary{decorations.markings}
\usetikzlibrary{arrows}
\usetikzlibrary{calc}
\usetikzlibrary{decorations.markings}

\tikzstyle{hvector}=[inner sep=2pt,draw=blue!50,fill=blue!10,thick]
\tikzstyle{unit}=[inner sep=2pt,shape=circle, draw]
\tikzstyle{counit}=[inner sep=2pt,shape=circle, draw,fill=gray]
\tikzstyle{antipode}=[inner sep=2pt,shape=rectangle, draw]
\tikzstyle{cocycle}=[inner sep=2pt,shape=circle, draw]
\tikzstyle{twistedm}=[inner sep=2pt,shape=circle, fill=gray]
\tikzstyle{autom}=[inner sep=2pt,shape=circle, draw]
\tikzstyle{coact}=[inner sep=2pt,shape=circle, fill=black]

\usepackage[bookmarks=true,%
    colorlinks=true,%
    linkcolor=blue,%
    citecolor=blue,%
    filecolor=blue,%
    menucolor=blue,%
    urlcolor=blue,%
    breaklinks=true]{hyperref}
\usepackage{slashed}    
\usepackage{verbatim}
\usepackage[normalem]{ulem}  
\usepackage{diagbox}
\usepackage[shortlabels]{enumitem}
\usepackage{tikz,tikz-cd}
\usetikzlibrary{knots, hobby, decorations.pathreplacing,
  shapes.geometric, calc}
\usepackage{orcidlink}

\newtheorem{theorem}{Theorem}[section]
\theoremstyle{definition}
\newtheorem{proposition}[theorem]{Proposition}

\newtheorem{definition}[theorem]{Definition}
\newtheorem{remark}[theorem]{Remark}

\newtheorem{question}[theorem]{Question}

\newtheorem{example}[theorem]{Example}

\def\BZ{\mathbbm Z}
\def\BQ{\mathbbm Q}

\def\BC{\mathbbm C}

\def\be{\begin{equation}}
\def\ee{\end{equation}}

\def\Aut{\mathrm{Aut}}
\def\LG{\mathrm{LG}}

\makeatletter

\renewcommand\thepart{\@Roman\c@part}%
\renewcommand\part{%
   \if@noskipsec \leavevmode \fi
   \par
   \addvspace{6.7ex}%
   \@afterindentfalse
   \secdef\@part\@spart}
\def\@part[#1]#2{%
    \ifnum \c@secnumdepth >\m@ne
      \refstepcounter{part}%
      \addcontentsline{toc}{part}{Part~\thepart.\ #1}%
    \else
      \addcontentsline{toc}{part}{#1}%
    \fi
    {\parindent \z@ \raggedright
     \interlinepenalty \@M
     \normalfont
     \ifnum \c@secnumdepth >\m@ne
       \centering\large\scshape \partname~\thepart.%
       \hspace{1ex}%
     \fi%
     \large\scshape #2%
     \markboth{}{}\par}%
    \nobreak
    \vskip 4.7ex
    \@afterheading}
  \def\@spart#1{
  \refstepcounter{part}%
  \addcontentsline{toc}{part}{#1}%
    {\parindent \z@ \raggedright
     \interlinepenalty \@M
     \normalfont
     \centering\large\scshape #1\par}%
     \nobreak
     \vskip 4.7ex
     \@afterheading}
\renewcommand*\l@part[2]{%
  \ifnum \c@tocdepth >-2\relax
    \addpenalty\@secpenalty
    \addvspace{0.75em \@plus\p@}%
    \begingroup
      \parindent \z@ \rightskip \@pnumwidth
      \parfillskip -\@pnumwidth
      {\leavevmode
       \normalsize \bfseries #1\hfil \hb@xt@\@pnumwidth{\hss #2}}\par
       \nobreak
       \if@compatibility
         \global\@nobreaktrue
         \everypar{\global\@nobreakfalse\everypar{}}%
      \fi
    \endgroup
  \fi}

\def\l@subsection{\@tocline{2}{0pt}{2pc}{6pc}{}}
\makeatother

\begin{document}

\title[Patterns of the $V_2$-polynomial of knots]{
Patterns of the $V_2$-polynomial of knots}
\author{Stavros Garoufalidis}
\address[Stavros Garoufalidis]{
  International Center for Mathematics, Department of Mathematics \\
  Southern University of Science and Technology \\
  Shenzhen, China \newline
  {\tt \url{http://people.mpim-bonn.mpg.de/stavros}} }
\email{stavros@mpim-bonn.mpg.de}
\author{Shana Yunsheng Li}
\address[Shana Yunsheng Li]{
  Department of Mathematics \\
  University of Illinois \\
  Urbana, IL, USA \newline
  {\tt \url{https://li-yunsheng.github.io}}}
\email{yl202@illinois.edu}

\date{16 March 2026}

\begin{abstract}
  Recently, Kashaev and the first author constructed an $R$-matrix from a
  Nichols algebra with an automorphism, that leads, via the Reshetikhin--Turaev
  functor, to a multivariable polynomial invariant of knots. Applying this to a
  rank 2 Nichols algebra, results in a sequence $V_n$ of 2-variable
  knot polynomials with integer coefficients, the first polynomial been
  identified with the Links--Gould
  polynomial. In this note we present the results of the computation of the
  $V_n$-polynomials for $n=1,2,3,4$. This leads to the discovery of emerging patterns,
  including the genus bound for $V_2$ being an equality for all 352.2 million 
  knots with at most $19$ crossings, as well as 
  unexpected Conway mutations that seem undetected by the $V_n$-polynomials
  as well as by Heegaard Floer Homology and Khovanov Homology.
\end{abstract}

\maketitle

{\footnotesize
\tableofcontents
}


\section{Introduction}
\label{sec.intro}

\subsection{A sequence of 2-variable knot polynomials}
\label{sub.intro}

Recently, Kashaev and the first author constructed an $R$-matrix from a
Nichols algebra with an automorphism~\cite{GK:multi}, that leads,
via the Reshetikhin--Turaev functor~\cite{RT:ribbon}, to a multivariable polynomial
invariant of knots. The definition of the $R$-matrix associated to a Nichols algebra
with automorphism (or its modules) is explicit~\cite[Thm.3.5,3.6]{GK:multi},
and uses the structure constants of the multiplication, comultiplication, evaluation,
coevaluation, antipode and automorphism of the Nichols algebra with automorphism.
Since Nichols algebras are not widely known in quantum topology, a reference for
them is Andruskiewitsch--Schneider~\cite{AS:pointed}.

Applying the construction of~\cite[Sec.7]{GK:multi} to one of
the simplest rank 2 Nichols algebra of diagonal type with automorphism, results
in a sequence of explicit $R$-matrices (on a $4(n+1)$-dimensional vector space)
for $n \geq 1$ and hence to a sequence $V_n(t,q)$ of 2-variable knot polynomials
with integer coefficients.
This algebra depends on one variable $q$ that determines the braiding and two
variables $(t_1,t_2)$ that determine the automorphism type. When $q$ is not
a root of unity and $(t_1,t_2,q)$ satisfy the relation $t_1t_2q^n=1$ for some
positive integer $n>0$, then the Nichols algebra has a right Yetter--Drinfel'd
f-module $Y_n$ of dimension $4n$~\cite[Prop.7.4]{GK:multi} and an explicit
$R$-matrix $T_n$.

Taking this as a black box, and parametrizing the three variables
$(t_1,t_2,q)$ satisfying $t_1t_2q^n=1$ in terms of two variables $(t,q)$ via 
$(t_1,t_2)=(1/(q^{n/2} t),t/q^{n/2})$, it was shown in~\cite{GK:multi} that $Y_n$
comes equipped with an $R$-matrix $T_n$ which leads to a matrix-valued knot invariant
$K \mapsto J_{T_n}(K) \in \operatorname{End}(Y_n)$
as well as to a scalar-valued invariant
$V_{K,n}(t,q)\in\BZ[t^{\pm 1},q^{\pm 1/2}]$ given by the $(1,1)$-entry of
$J_{T_n}(K)$.

It was advocated in~\cite{GK:multi} that the sequence $V_n$ of 2-variable knot
polynomials has similarities and differences with the sequence of the Jones
polynomial of a knot and its parallels, otherwise known as the colored Jones
polynomial. Similarities include their specialization to the Alexander polynomial
discussed shortly (which is part of an MMR Conjecture and of a loop expansion) and
$q$-holonomicity, whereas the main difference is that the $V_n$-polynomials involve
a variable $t$ lacking in the colored Jones polynomials, and it is the degree with
respect to this variable that is related to the genus bounds discussed below.

The polynomial invariant $V_n$ satisfies the symmetry~\cite{GK:multi}
\be
\label{Vnsym}
V_{K,n}(t,q)=V_{K,n}(t^{-1},q), \qquad V_{\overline K,n}(t,q)=V_{K,n}(t,q^{-1}) 
\ee
where the first equality comes from the involution exchanging $t_1$ and $t_2$,
and in the second one $\overline K$ denotes the mirror image of $K$. Three further
properties of the $V_n$-polynomials were conjectured in~\cite{GK:multi}, namely

\noindent
$\bullet$ 
The specialization
\be
\label{Vnspecial}  
V_{K,n}(q^{n/2},q)=1, \qquad V_{K,n}(t,1) =\Delta_K(t)^2
\ee
where $\Delta_K(t) \in \BZ[t^{\pm 1/2}]$ is the symmetrized
(i.e., $\Delta_K(t)=\Delta_K(t^{-1})$), normalized (i.e., $\Delta_K(1)=1$)
Alexander polynomial.

\noindent
$\bullet$ 
The relation
\be
\label{VLG1}
V_1=\LG
\ee
with the Links--Gould polynomial invariant~\cite{Links-Gould}.

\noindent
$\bullet$ 
The genus bound
\be
\label{degVn}
\deg_t V_{K,n} \leq 4 g(K)
\ee
where the Seifert genus $ g(K)$ is the smallest genus of a spanning surface of a knot.
Here, by $t$-degree $\deg_t$ of a Laurent polynomial of $t$ we mean the difference
between the highest and the lowest power of $t$.
We will say that $V_n$ detects the genus of a knot $K$ if~\eqref{degVn} is an equality.

The paper~\cite{GK:multi} stimulated a lot of subsequent recent work. 
The relation $V_1=\LG$ is now known and it follows from the fact
that both polynomials satisfy a cubic skein theory with a unique solution with
value $1$ at the unknot and vanishing for split links~\cite{LG-V1}. The relation
~\eqref{VLG1}, combined with the work of Kohli--Tahar~\cite{Kohli-Tahar} who proved
a genus bound for the Links--Gould polynomial, implies the genus bound~\eqref{degVn}
for $n=1$. 

What's more, the relation~\eqref{VLG1} lead to a conjecture
\be
\label{VLGn}
V_n = \LG^{(n)}
\ee
that relates the $V_n$-polynomials of knots with the colored LG-polynomials
$\LG^{(n)}$ of knots, where the latter are the $\mathfrak{sl}(2|1)$-quantum group
invariants of knots colored by a simple $4n$-dimensional
representation~\cite[Conj.1.4]{Vn}.

The above conjecture~\eqref{VLGn} identifies the Nichols algebra polynomial
invariants $V_n$ with the LG-polynomials of a knot and its parallels
(see~\cite[Thm.1.1]{Vn}), and explains that the $V_n$-polynomials are Vassiliev
power series invariants. What's more, the specialization and the genus bound
properties for the $\LG^{(n)}$ polynomials are known (see~\cite[Thm.1.3]{Vn}),
hence Conjecture~\eqref{VLGn} implies~\eqref{Vnspecial} and~\eqref{degVn} for all
$n \geq 1$.

Finally, Conjecture~\eqref{VLGn} is known for $n=1$
(as mentioned above~\cite[Thm.1.1]{LG-V1}) and also for $n=2$ (see~\cite[Thm.1.5]{Vn}).

With regard to the genus bound~\eqref{degVn} for $n=1$, 
$V_1=\LG$ cannot detect the genus of a knot since it is known that the Links--Gould
invariant does not detect Conway mutation, whereas the genus does. On the other hand,
it follows from~\cite{Vn} that $V_2$ can be expressed in terms of $V_1$ polynomial
of a knot and its $(2,1)$-parallel, and this leads to a proof of the
genus bound~\eqref{VLGn} for $n=2$.

Of all those properties of $V_2$, the most intriguing one is whether it detects the
genus of all knots. Based on some experiments with a few knots of 12 and 13 crossings,
it was observed in~\cite{GK:multi} that in all computed cases $V_2$ detects the genus
of the sampled knots. Is this an accident for knots with low numbers of
crossings, or a new phenomenon? To decide one way or another, one needs an efficient
way to compute the $V_n$-polynomials of knots, do so, and sieve the data.
This is exactly what we did, and led to the results of our paper. We remark that
the method we used here can be easily promoted to a subexponential fixed-parameter
tractable algorithm for computing knot polynomials coming from the
Reshetikhin--Turaev functor \cite{Shana:FPT-RT}.

In the next sections we present the results of our computations and the patterns found,
and in Section~\ref{sec.compute} we present the details of the algorithm. All data
are available in~\cite{GS:VnData}.

\subsection{Does $V_2$ detect the genus of a knot?}
\label{sub.results}


Since we are talking about tables of knots and their invariants, we will be 
using the naming of the HTW table of knots up to 16 crossings~\cite{HTW}
imported in \texttt{SnapPy}~\cite{snappy} and also in
\texttt{KnotAtlas}~\cite{knotatlas}. In addition, we used Burton's table and
notation \cite{Burton:knots} for the list of knots with $17$, $18$ and $19$ crossings. 
The findings are summarized in Table~\ref{tab:1} below.



\begin{small}
\begin{table}[htpb!]
\begin{center}
\begin{tabular}{|c|r|r||r|r|r|}
\hline
\begin{small} crossings \end{small} & 
\begin{small} Knots \end{small} &
\begin{small} $\Delta=1$ \end{small} &
\begin{small} $\Delta$ fails  \end{small} &
\begin{small} $V_1$ fails \end{small} &                                        
\begin{small} $V_2$ fails \end{small}  \\ \hline
$\leq 10$ & 249 & 0 & 0 & 0 & 0
\\ 
11 & 552 & 2 & 7 & 7 & 0
\\ 
12 & \num{2176} & 2 & 29 & 20 & 0
\\ 
13 & \num{9988} & 15 & 208 & 173 & 0
\\ 
14 & \num{46972} & 36 & \num{1220} & 974 & 0
\\ 
15 & \num{253293} & 118 & \num{6319} & \num{5025} & 0
\\ 
16 & \num{1388705} & 499 & \num{48174} & \num{37205} & 0
\\ 
17 & \num{8053393} & \num{1734} & \num{303823} & \num{228996} & 0
\\ 
18 & \num{48266466} & \num{6850} & \num{2001954} & \num{1481428} & 0
\\ 
19 & \num{294130458} & \num{25647} & \num{13287958} & \num{9676780} & 0 \\
\hline
\end{tabular}
\end{center}
\caption{Knot counts, up to mirror image. The three ``fail'' columns list the number
  of knots where the corresponding polynomials fail to detect their genus.}
\label{tab:1}
\end{table}
\end{small}

The inequality~\eqref{degVn} combined with the specialization~\eqref{Vnspecial}
for $n=2$ imply that
\be
2 \deg_t \Delta_K(t) \leq \deg_t V_{K,2}(t,q) \leq 4 g(K) \,.
\ee
On the other hand,
\be
\deg_t \Delta_K(t) \leq 2 g(K)
\ee
with equality if and only if a knot is Alexander-tight, otherwise Alexander-loose.
In our paper we abbreviate these two classes simply with tight/loose, similar to what
people do in Heegaard Floer Homology (abbreviated by HFK) \cite{SS:thick} and Khovanov Homology
where they talk about
HFK-thin/thick or Kh-thin/thick knots, but then they drop the HFK or Kh once
the context is clear. Likewise, we use the term thin in our paper to mean HFK and
Kh-thin. Note that alternating knots are tight~\cite{Murasugi}.
What's more, quasi-alternating knots (a class that includes all alternating knots)
introduced by Ozsv\'ath--Szab\'o in~\cite{OS:quasi} are HFK and
Khovanov-thin~\cite{Manolescu:thin}, and hence tight.

Combining the above two inequalities, it follows that
the inequality~\eqref{degVn} is in fact an equality for $V_2$ and all tight knots.
Note next that there are no loose knots with $\leq 10$ crossings. Moreover, the number
of loose knots with $\leq 16$ crossings is given in Table~\eqref{tab:1}.
Incidentally, the list of loose knots was compiled by computing in \texttt{SnapPy}
the Alexander polynomial, and also the HFK (and in particular, the Seifert genus
of a knot).

Among the loose knots, are the ones with trivial Alexander polynomial (also computed
by \texttt{SnapPy}) which are in some sense extreme. The list of 7 loose knots
(up to mirroring) with 11 crossings is
\be
\label{loose11}
11n34^\ast, \hspace{0.2cm} 11n42^\ast, \hspace{0.2cm} 11n45, \hspace{0.2cm} 11n67,
\hspace{0.2cm} 11n73, \hspace{0.2cm} 11n97, \hspace{0.2cm} 11n152 
\ee
where the asterisque indicates that the knot has trivial Alexander polynomial, and
the pair $(11n34,11n42)$ is the famous Kinoshita--Terasaka and Conway pair of
mutant knots. Their genuses are $3, 2, 3, 2, 3, 2, 3$ respectively, $t$-degrees of the
$V_1$-polynomial are $6, 6, 8, 6, 8, 6, 8$ and $t$-degrees of the $V_2$-polynomial
are $12, 8, 12, 8, 12, 8, 12$, confirming the equality in~\eqref{degVn} for $n=2$. 

\begin{table}[htpb!]
\begin{center}
\begin{tabular}{|c|r|r|r|r|}
\hline    
polynomial & $V_1$ & $V_2$ & $V_3$ & $V_4$ \\
\hline
Knots & $\leq 16$ & $\leq 16$ & $\leq 14$ & $\leq 14$  \\
\hline
Loose knots &  $\leq 19$ & $\leq 17$ &  &  \\
\hline  
\end{tabular}
\end{center}
\caption{Table of computed values of $V_n$, the values indicate the number of
  crossings.}
\label{tab:2}
\end{table}

Table~\ref{tab:2} summarizes the knots for which the $V_n$-polynomial was computed.
The data is available in~\cite{GS:VnData}, with the convention that we replaced
$q$ by $q^2$ in $V_1$ and $V_3$ so that we obtain Laurent polynomials in $t$ and $q$,
as opposed to Laurent polynomials in $t$ and $q^{1/2}$. 
For $n=2$, we computed its values for all loose knots with at most 17 crossings.
In all cases, we found that the inequality~\eqref{degVn} is an equality for $n=2$.
Combined with the specialization and the genus bounds for $V_2$, this implies
the following.

\begin{proposition}
\label{prop.1}
$V_2$ detects the genus of all $\num{352152252}$ knots with at most $19$ crossings.
\end{proposition}

We remark that the results for 18 and 19 crossing knots involve extra structure of the
$V_n$-polynomials beyond the 1-loop specialization of Equation~\eqref{Vnspecial}.
We will explain this extra structure concerning the 2-loop specialization of the
$V_n$-polynomials in a subsequent publication.

\begin{remark}
\label{rem.genus}  
The relation between $V_2$ and $V_1$ discussed in Section~\ref{sub.relV1V2} below,
combined with the fact that $V_1=\LG$ implies that the map
\be
\label{Vas}
K \mapsto V_{K,2}(e^{\hbar N},e^\hbar) \in \BQ[N][\![\hbar]\!]
\ee
is a Vassiliev power series invariant of knots. Hence, if ~\eqref{degVn} is an
equality for $n=2$, it follows that Vassiliev invariants determine the Seifert genus
of a knot, and hence (answering a folk conjecture) detect the unknot. 
A celebrated method to detect the genus of the knot is Heegaard Floer
Homology~\cite{OS}. A second (conjectural) method to compute the genus of a knot uses
hyperbolic geometry, and more specifically the degree of the twisted torsion
polynomial $\tau_{K,2}(t)$ of a hyperbolic knot (twisted with the
adjoint representation of the geometric representation of a hyperbolic knot); see
~\cite[Sec.1.6]{Dunfield:twisted}. Note that Conjecture 1.7
of~\cite{Dunfield:twisted} was verified for all hyperbolic knots with
at most 15 crossings. 
\end{remark}

\subsection{$V_2$ detects the genus for the Kinoshita--Terasaka family}
\label{sub.KT}

In this section we study the genus bound for $V_2$
for a 2-parameter family $(KT_{r,n},C_{r,n})$ of pairs of mutant knots with trivial
Alexander polynomial, with Seifert genus $g(KT_{r,n})=r$ and $g(C_{r,n})=2r-1$.
The family was introduced by Kinoshita--Terasaka~\cite{KT} and further studied
by Ozsv\'ath--Szab\'o~\cite{OS:genusbounds} as an application of the Heegaard Floer
Homology of knots. The computation of the genus is due to Gabai~\cite{Gabai},
and the family
is discussed in detail in~\cite{OS:genusbounds}, where a drawing may be found in
Fig.1 of~\cite{OS:genusbounds}. The planar projection of each of the knots
$KT_{r,n}$ and $C_{r,n}$ requires $4r+4n+2$ crossings, and the case of $(r,n)=(2,1)$
is the famous Kinoshita--Terasaka and Conway pair of knots.

We computed the $V_2$ polynomial for $r=2,\dots,20$ and $n=2$ and its $t$-degree
exactly matches the genus. 

%

Since these knots involve a full $n$-twist on two strands, general TQFT reasons
imply that their $V_2$-polynomial, as a function of $n$, satisfies a linear
recursion relation with coefficients the eigenvalues of the square of the $R$-matrix.
In our case, there are 7 eigenvalues.
This and the computation of the $V_2$ polynomial
for a fixed value of $r$ and 7 consecutive values of $n$ gives the $V_2$-polynomial
for fixed $r$ and all $n$. Computing so for $r=1,\dots,10$ and $n=1,\dots,8$, we
verified that the degree of the $V_2$ polynomial matches the genus for 
$r=1,\dots,10$ and $n \in \BZ$.

Summarizing, we obtained the following.

\begin{proposition}
\label{prop.KT}
Equality holds in~\eqref{degVn} for all pairs of knots $(KT_{r,n},C_{r,n})$
for $r=2,\dots,20$ and $n=2$, as well as for $r=1,\dots,10$ and $n \in \BZ$.
\end{proposition}
  
\subsection{$V_2$-trivial Conway mutations}
\label{sub.Vtrivial}

A question that we discuss next is how strong is the new $V_2$ polynomial in separating
knots. Given the values of the polynomial for knots up to 14 crossings, we searched
for repetitions, taking into account mirror image, which changes $V_2(t,q)$ to
$V_2(t,q^{-1})$. Here, we came across a new surprise. The $V_2$ polynomial separates
knots with at most 11 crossings, but fails to separate 12 crossing knots, and
the three pairs that we found are
\be
\label{pairs12a}
\begin{aligned}
(12n364, \,\,\, \overline{12n365})  &&
(12n421, \,\,\, \overline{12n422})  &&
(12n553, \,\,\, \overline{12n556}) \,. 
\end{aligned}
\ee
We tried the $V_1$ polynomial on them and it failed to
separate them, and we tried the $V_3$-polynomial which also failed to separate. Yet,
the genus inequality~\eqref{degVn} was an equality for $n=2$, which meant that these
3 pairs have equal genus (in each pair). We checked their HFK,
computed by \texttt{SnapPy}, and their Khovanov Homology, computed by
\texttt{KnotAtlas}, which, a bit to our surprise, were equal in each pair.
Looking at these 3 pairs more closely, we realized that they are in fact
Conway mutants. A table of mutant knots with at most 15 crossings is given in
Stoimenow~\cite{Stoimenow}. As was pointed out to us by N. Dunfield, one
can separate the knots in these pairs using the homology of their 5-fold covers,
or the certified isometry signature of the complete hyperbolic structure.

Having found these unexpected pairs of Conway mutant knots, we tried knots with
13 crossings, where we now found 50 more pairs with exactly the same properties as
above, given in Table~\eqref{pairs13}. We then searched knots with
14 crossings, where we now found $333$ pairs and $1$ triple (up to mirror image)
with the same properties as above.

For lack of a better name, let us say that two knots are $V_2$-equivalent
if they have equal $V_2$-polynomial; this notion is similar to the almost-mutant
knots of~\cite{Dunfield:mutation}. The $V_2$-equivalence classes of knots of size
more than 1 with at most 12, 13 and 14 crossings are given in Tables~\eqref{pairs12},
~\eqref{pairs13} and~\eqref{pairs14a}, respectively. A computer readable list
of the $V_2$-equivalent tuples, up to $16$ crossings,
is the file \texttt{V2-equiv-tuples.txt} available at \cite{GS:VnData}. The numbers of $V_2$-equivalence tuples are summarized in the following table.

\begin{table}[htpb!]
\begin{center}
\begin{tabular}{|r|r|r|r|r|r|r|}
\hline    
crossings & $\leq 11$ & 12 & 13 & 14 & 15 & 16 \\
\hline
& $\emptyset$ & $2\!:\!3$ & $2\!:\!50$ & $2\!:\!333, \, 3\!:\!1$
& $2\!:\!2324, \, 3\!:\!38$ & $2\!:\!14387, \, 3\!:\!214, \, 4\!:\!17, \, 6\!:\!8$ \\
\hline  
\end{tabular}
\end{center}
\caption{Number of $V_2$-equivalence classes of size more than 1 (up to mirror
  image) where $n\!:\!d$ means the number of $n$-tuples is $d$.}
\label{tab:3}
\end{table}


Most of the $V_2$-equivalent knots that we found have equal $V_n$-polynomials for
$n=1,2,3,4$. Let us say that two knots are $V_2$-exceptional if they have equal
$V_1$ and $V_2$ polynomials but different $V_3$ or $V_4$ polynomials.

\begin{table}[htpb!]
\tiny
\begin{equation*}
\begin{aligned}
  (14n2420,14n4659)^3 && (14n2423, 14n5868)^2 && (14n5822, 14n5852)^3 &&
  (14n5828,14n6370)^2 \\
  (15n11298, \,15n21931)^3 &&
(15n11303, \,15n29554)^2 &&
(15n29403, \,15n29499)^3 &&
(15n29411, \,15n33147)^2 \\
(15n107431, \,15n107988)^3 &&
(15n107943, \,15n109145)^3 && 
(16n2191, \,16n24753)^2 &&
(16n2218, \,16n24757)^2 \\
(16n2670, \,16n24603)^2 &&
(16n2675, \,16n24611)^2 &&
(16n4803, \,16n24709)^3 &&
(16n4809, \,16n24717)^3 \\
(16n4829, \,16n24799)^3 &&
(16n4833, \,16n24803)^3 &&
(16n58598, \,16n112012)^3 &&
(16n58607, \,16n154441)^2 \\
(16n58979, \,16n112235)^3 &&
(16n58982, \,16n112242)^3 &&
(16n59065, \,16n179415)^2 &&
(16n59066, \,16n179416)^2 \\
(16n72711, \,16n90534)^3 &&
(16n72714, \,16n90535)^3 &&
(16n112158, \,16n154579)^3 &&
(16n112160, \,16n154580)^3 \\
(16n112166, \,16n179120)^2 &&
(16n112167, \,16n179121)^2 &&
(16n153640, \,16n154120)^3 &&
(16n153647, \,16n178690)^2 \\
(16n178206, \,16n178434)^3 &&
(16n178208, \,16n178442)^3 &&
(16n178307, \,16n191518)^2 &&
(16n178308, \,16n191519)^2 \\
(16n178384, \,16n178847)^3 &&
(16n178386, \,16n178850)^3 &&
(16n178397, \,16n190704)^2 &&
(16n178398, \,16n190705)^2 \\
(16n272487, \,16n357198)^3 &&
(16n323846, \,16n376154)^3 &&
(16n396926, \,16n400281)^2 &&
(16n401578, \,16n402067)^2 \\
(16n657584, \,16n661016)^3 &&
(16n666167, \,16n666234)^3 &&
(16n863741, \,16n906884)^2 &&
(16n906885, \,16n906887)^2 
\end{aligned}
\end{equation*}
\normalsize
\caption{$V_2$-exceptional pairs with equal $V_1$ and $V_2$ polynomials but different
  $V_3$ and $V_4$ polynomials. The superscripts refer to the flavors of the pairs,
  as defined in \eqref{3flavors}.}
\label{tab:exceptional}
\end{table}

Summarizing, we obtain the following.

\begin{proposition}
\label{prop.2a}
Up to $16$ crossings, $V_2$-equivalent knots have equal $V_1$, $V_2$, $V_3$ and
$V_4$-polynomials, except those in Table~\ref{tab:exceptional}.
\end{proposition}

\begin{proposition}
\label{prop.2b}
Up to $15$ crossings, $V_2$-equivalent knots are Conway mutant. Hence 
they have equal colored Jones polynomials, ADO polynomials and HOMFLY polynomial,
and when hyperbolic, they have equal volume and trace field.
\end{proposition}

The above proposition is expected to hold for up to $16$ crossings, only that the
table of Conway mutant 16 crossing knots is not known. In particular, we have
verified that all $V_2$-equivalent knots, up to $16$ crossings, have the same
HOMFLY polynomial.

We next checked the Khovanov Homology and the Knot Floer Homology of the
$V_2$-equivalent knots. Khovanov Homology in various versions, such as the original
one, the reduced, with coefficients in a field, is either known to be Conway mutation
invariant or expected to be so. For an updated discussion, see~\cite[Sec.9]{KWZ}
and references therein. In view of this, the next proposition was not a surprise.

\begin{proposition}
\label{prop.2c}
Up to $16$ crossings, $V_2$-equivalent knots have equal Khovanov Homology.
\end{proposition}


On the other hand, Conway mutation changes the genus of a knot, hence in general
changes HFK (which among other things, determines the genus of a knot), the
classic example being the Kinoshita-Terasaka and Conway pair of 11-crossing mutant
knots with trivial Alexander polynomial. Since $V_2$-equivalent knots appear to be
mutant, and mutation can change the genus, the genus equality~\eqref{degVn} may fail
for $n=2$. But this is not what happens in our sample.

\begin{proposition}
\label{prop.2d}
Up to $16$ crossings, $V_2$-equivalent knots have equal HFK.
\end{proposition}

We also remark that the two mutant knots $11n76$ and $\overline{11n78}$ have
equal HFK, but different $V_2$-polynomials. Hence $V_2$ is not determined
by HFK and detects some mutations that HFK misses.

\begin{remark}
Recall that the canonical genus of $K$ is the minimal genus of all Seifert surfaces
produced by applying Seifert's algorithm on all diagrams of $K$. The $z$-degree of
the HOMFLY polynomial $P_K(a,z)$, divided by $2$, gives a lower bound $L_K$ of the
canonical genus of knot $K$ \cite{Morton:cg}. When $L_K$ is strictly larger than the
genus of $K$, the genus of $K$ cannot be realized by Seifert surfaces given by
Seifert's algorithm. Among the $17375$ $V_2$-equivalence classes up to $16$ crossings,
only $83$ of them have $L_K$ strictly larger than their genuses, and all $83$ of them
are thick. 
  
On the other hand, all $48$ $V_2$-exceptional pairs in Table~\ref{tab:exceptional}
are thick, with $22$ being tight and $26$ being loose. Taking a closer look, it turns
out that the $22$ tight and thick $V_2$-exceptional pairs form exactly the
intersection of the $48$ $V_2$-exceptional pairs and the $83$ $V_2$-equivalence
classes with $L_K$ strictly larger than genus. 
  
In short, among all $17375$ $V_2$-equivalence classes up to $16$ crossings, 
\begin{equation}
\begin{aligned}
\{V_2\text{-exceptional}\} \cup \{L_K > \text{ genus}\} \subset \{\text{thick}\},\\
\{V_2\text{-exceptional}\}\setminus \left( \{\text{tight}\}\cap\{\text{thick} \}\right)\neq \varnothing,\\
\{L_K > \text{ genus}\}\setminus \left(\{\text{tight}\}\cap\{\text{thick} \}\right) \neq \varnothing,
\end{aligned}
\end{equation}
while
\begin{equation}
  \{V_2\text{-exceptional}\}\cap \{\text{tight}\}\cap\{\text{thick} \} =
  \{V_2\text{-exceptional}\}\cap \{L_K > \text{ genus}\} \neq \varnothing.
\end{equation}
Given the substantial size of the $V_2$-equivalence classes in consideration, this
is not likely to be random.
\end{remark}

\begin{question}
\label{que.1}
Are $V_2$-equivalent knots always Conway mutant? Do they always have
equal HFK and Khovanov Homology? And why?
\end{question}

We can give a partial answer to this question as follows. The observed
$V_2$-equivalence classes come in 3 flavors

\be
\label{3flavors}
1: \text{tight + thin}, \qquad 2: \text{tight + thick},
\qquad 3: \text{loose + thick} 
\ee
and the counts of the $V_2$-equivalence classes of size more than one,
of knots with $\leq 14$ crossings, according to their flavor is given
by 
\be
\label{flavor.count}
1 : 293, \qquad 2 : 57, \qquad 3 : 37.
\ee


Regarding the more numerous class 1, note that the HFK homology of an HFK-thin
(resp., Kh-thin) knot is determined by its Alexander polynomial (resp., by the
Jones polynomial and the signature). Since mutation does not change the Alexander
polynomial, nor the Jones polynomial, nor the signature, it follows that mutant
thin knots have equal HFK and equal Khovanov Homology. This gives an explanation
of the last two parts in Question~\ref{que.1} for the class 1, which as was
mentioned in the introduction, includes pairs of mutant quasi-alternating knots.

The tuples in the other two classes are not as well-understood. Our tables
(given in the Appendix) give concrete examples of tight + thick or loose + thick knots.
For instance, three pairs of tight + thick knots are
\be
\label{flavor2.ex}
(13n1655, \,\,13n1685)^2, \qquad
(14n1370, \,14n1395)^2, \qquad
(14n1699, \,14n1947)^2
\ee
and three pairs of loose + thick knots are
\be
\label{flavor3.ex}
(13n372, \,\,13n375)^3, \qquad
(13n536, \,\,13n551)^3, \qquad
(13n1653, \,\,13n1683)^3 \,.
\ee
There are several methods of constructing knots with equal HFK and equal Khovanov
Homology discussed for example in detail in Hedden--Watson,~\cite{Hedden:geography},
but we do not know how to apply these constructions to generate our examples.

It is also worth noticing that all tight + thin knots listed in the appendix,
except $14n14135$, are $\overline{\text{Kh}'}$-thin, as defined in
\cite[Def.5.1]{ORS:ODD}.
\begin{figure}[htpb!]
\centering
\scalebox{.7}{
$\vcenter{\hbox{\includegraphics{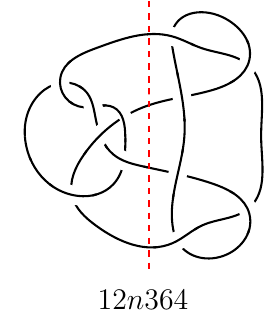}}}$

\qquad
$\vcenter{\hbox{\includegraphics{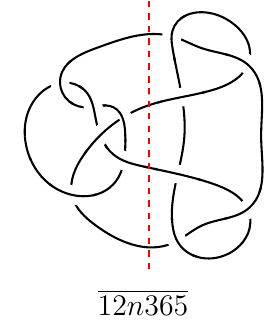}}}$

\qquad

$\vcenter{\hbox{\includegraphics{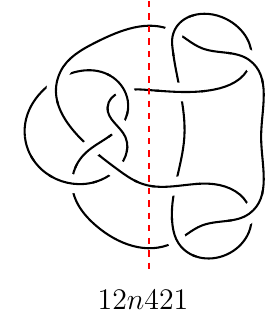}}}$
\qquad
$\vcenter{\hbox{\includegraphics{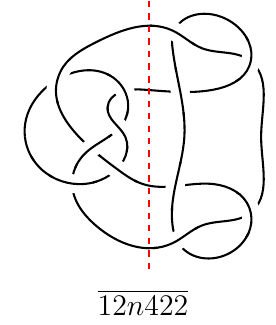}}}$
}
\\

\scalebox{.7}{
$\vcenter{\hbox{\includegraphics{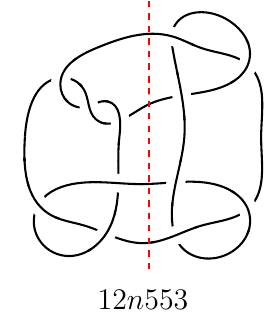}}}$
\qquad
$\vcenter{\hbox{\includegraphics{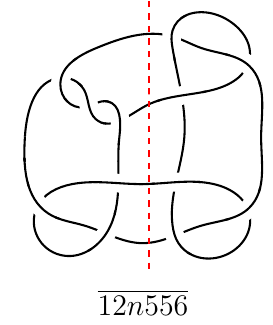}}}$
}

\caption{The 3 pairs of knots from~\eqref{pairs12a}.}
\label{f.3knots}
\end{figure}

\subsection{Independence of the $V_2$ and the 2-loop polynomials}
\label{sub.2loop}

The colored Jones polynomial can be decomposed into loop invariants, starting
from $0$-loop which is the inverse Alexander polynomial, and then going up the loops.
In fact, the 2-loop invariant of the Kontsevich integral of a knot is essentially a
2-variable polynomial invariant, and its image under the $\mathfrak{sl}_2(\BC)$-weight
system is a 1-variable polynomial computed efficiently by Bar-Natan and van der Veen
~\cite{BNV:PG,BNV:API}.

One may ask whether the $V_2$-polynomial is determined by the 2-loop invariant $Z_2$
of the Kontsevich integral. It is known that the map $K \mapsto Z_2(\mathrm{Wh}(K))$
is a degree 2 Vassiliev invariant of knots~\cite[Thm.1]{Ga:whitehead}, where
$\mathrm{Wh}(K)$ denotes the Whitehead doubling of a 0-framed knot with a positive
clasp. Since the vector space of degree 2 Vassiliev invariants is 1-dimensional
generated by $a_2(K)=\Delta''(1)$, it follows that
$Z_2(\mathrm{Wh}(K)) = c \, a_2(K)$ for a universal nonzero constant $c$.
For the trefoil and its mirror image, we have $a_2(3_1)=a_2(\overline{3_1})=2$, which
implies that
\be
\label{Z2}
Z_2(\mathrm{Wh}(3_1)) = Z_2(\mathrm{Wh}(\overline{3_1})) \,.
\ee
On the other hand, $V_{\mathrm{Wh}(3_1),2}(t,q) \neq V_{\mathrm{Wh}(\overline{3_1}),2}(t,q)$,
the exact values given in Equation~\eqref{Wh31} in the Appendix. This implies the
following. 

\begin{proposition}
The $V_2$-polynomial is not determined by the 2-loop part of the Kontsevich integral.   
\end{proposition}

Based on limited computations available, it was observed in~\cite{GK:multi} that
Khovanov Homology alone, or HFK alone, or the colored Jones polynomial alone
do not determine $V_2$. 

\subsection{A relation between $V_1$ and $V_2$}
\label{sub.relV1V2}

In a sense, the sequence of $V_n$-polynomials are similar to the sequence of the
Jones polynomial of a knot~\cite{Jones} and its $(n,1)$-parallels. In fact, it
follows from the axioms of the TQFT and the tensor product decomposition of two
irreducible representations of $\mathfrak{sl}_2$ (known as the Clebsch-Gordan formula)
that the Jones polynomial of a parallel of
knot is a linear combination (with coefficients that are independent of the knot) of 
colored Jones polynomial, colored by the irreducible representations of
$\mathfrak{sl}_2(\BC)$; see~\cite{RT:ribbon,Tu:book}, and vice-versa.

It was recently conjectured in~\cite{Vn} that the $V_n$-polynomials of a knot
are also linear combinations of the $V_1$-polynomial of a knot and its $(n',1)$-cables
for $n' \leq n$, and a proof for $n=2$ was given there. We illustrate this relation
here, at the same time giving a consistency between the coefficients of the
relation computed by the spectral decomposition of $R$-matrices and by
representation theory in~\cite{Vn} with the computer-program that computes $V_1$
and $V_2$. The following relation holds for the unknot,
$3_1$, $4_1$, $6_1$, $6_2$, $6_3$, $7_7$, $8_3$, $8_4$ and their mirrors (in total,
14 knots) and was conjectured to hold for all knots (now proven in~\cite[Sec.4.2]{Vn})

\be
\label{relV1V2}
V_{K,2}(t^2,q^2) = c_{2,0}(t,q) V_{K(2,1),1}(t,q)
+ c_{2,-1}(t,q) V_{K,1}(t^2 q^{-1},q)
+ c_{2,1}(t,q) V_{K,1}(t^2 q,q) \,,
\ee
where
\begin{small}
\be
\label{3c}
c_{2,-1}(t,q) = \frac{t(t^2 q^2 -1)}{q(1+q^2)(t^2-1)}, \qquad
c_{2,1}(t,q) = \frac{t^2-q^2}{qt(1+q^2)(t^2-1)}, \qquad
c_{2,0}(t,q) = \frac{(q+t)(1+qt)}{(1+q^2)t} 
\ee
\end{small}
satisfy the symmetry $c_{2,1}(t,q) = c_{2,-1}(t^{-1},q)$,
$c_{2,0}(t,q)=c_{2,0}(t^{-1},q)$.


Some values of~\eqref{relV1V2} are given in the appendix.

\subsection{$V_n$ detects the genus of torus knots with two strings}

From the general setting, it follows that if $\beta$ and $\gamma$ are elements of
a braid group of a fixed number of strands and $K_n$ denotes the link obtained 
by the closure of $\beta^n \gamma$, then $K_n$ is a knot if $n$ lies in an arithmetic
progression and the sequence $V_1(K_n)(t,q)$ is holonomic and satisfies a linear
recursion relation with coefficients in $\BZ[t^{\pm 1},q^{\pm 1}]$ coming from the
minimal
polynomial of the square of the $R$-matrix. This can be computed explicitly and
leads to the answer. The above holds locally, if we replace a tangle $\gamma$
in a planar projection of a knot by $\beta^n \gamma$, and holds for any of
the polynomial invariants that we discuss in this paper.

We illustrate this giving a recursion relation of the values of 
$V_1$ and $V_2$ for $T(2,2b+1)$-torus knots for an integer $b$.
The minimal polynomial of the square of the $R$-matrix of $V_1$ is
\be
(-1 + x) (-t^2 + q^2 x) (-1 + q^2 t^2 x) =
-t^2 + (q^2 + t^2 + q^2 t^4) x -(q^2 + q^4 t^2 + q^2 t^4) x^2 + q^4 t^2 x^3 \,.
\ee
It follows that $f_b(t,q) = V_{T(2,2b+1),1}(t,q)$ satisfies the recursion relation
\be
-t^2 f_b(t,q) + (q^2 + t^2 + q^2 t^4) f_{b+1}(t,q)
-(q^2 + q^4 t^2 + q^2 t^4) f_{b+2}(t,q) + q^4 t^2 f_{b+3}(t,q) =0
\ee
for $b \in \BZ$ with initial conditions
\be
\label{finit}
f_{-1}(t,q)=1, \quad f_0(t,q)=1, \quad f_{1}(t,q) =
1 + (q^{-1}+q^{-3}) u + q^2 u^2
\ee
where
\be
\label{u}
u=t+t^{-1}-q-q^{-1} \,.
\ee


This and the $ t \leftrightarrow t^{-1}$ symmetry of $V_1$ implies that
$f_b(t,q) = q^{-2b} (t^{2b}+t^{-2b}) + \text{(lower order terms)}$
for $b \geq 0$, thus $\deg_t(f_b(t,q))=4b=4 \cdot \text{genus}(T(2,2b+1))$ for
$b>0$. It follows that inequality~\eqref{degVn} for $n=2$ is an equality
for $b \geq 0$. Since $\overline{T(2,2b+1)}=T(2,-2b-1)$, it follows that
$f_b(t,q^{-1})=f_{-b-1}(t,q)$ which then concludes that inequality~\eqref{degVn}
for $n=1$ is an equality for all 2-string torus knots.

Likewise, the minimal polynomial of the square of the $R$-matrix of $V_2$ is
\be
(-1 + x) (-t^2 + q^2 x) (-1 + q^3 x) (-t^2 + q^4 x) (-1 + q^2 t^2 x) (-1 + q^4 t^2 x)
\ee
which translates into a 6th order linear recursion relation for 
$g_b(t,q) = V_{T(2,2b+1),2}(t,q)$ with initial conditions 

\be
\begin{aligned}
g_{-1}(t,q) = & g_0(t,q^{-1}) = 1 \\
g_{-2}(t,q) = & g_1(t,q^{-1}) =
1 + (q + 2 q^3 - q^4 + q^5 - q^6) u + (q^2 + q^4 - q^5) u^2
\\
g_{-3}(t,q) = & g_2(t,q^{-1}) =
1 + (2 q + 3 q^3 - q^4 + 3 q^5 - q^6 + 2 q^7 - q^8 + q^9 - 2 q^{10} + 
q^{11} - q^{12}) u \\
& + (4 q^2 + 7 q^4 - 3 q^5 + 10 q^6 - 6 q^7 + 
6 q^8 - 7 q^9 + 3 q^{10} - 3 q^{11}) u^2 \\
& + (3 q^3 + 6 q^5 - 3 q^6 + 6 q^7 - 6 q^8 + 3 q^9 - 3 q^{10}) u^3
+ (q^4 + q^6 - q^7 + q^8 - q^9) u^4
\end{aligned}
\ee
with $u$ as in~\eqref{u}. As in the case of $V_1$, from the above recursion
one deduces that inequality~\eqref{degVn} for $n=2$ is in fact an equality
for all 2-string torus knots. 

We have performed the analogous calculation for the case of the $V_3$ and $V_4$
polynomials, and the conclusion is that inequality~\eqref{degVn} for $n=1, \dots, 4$
is in fact an equality for all 2-string torus knots. Therefore

\begin{proposition}
  For $n=1,\dots, 4$, the $V_n$-polynomial detects the genus of all $2$-string torus knots.
\end{proposition}

\subsection{Positivity of the $V_1$ and $V_2$-polynomials of
alternating knots?}
\label{sub.positivity}

The next topic that we discuss is a curious positivity observation for the
coefficients of the $V_1$ and $V_2$ polynomials of alternating knots.
Recall the number of alternating knots with at most 16 crossings
(up to mirror image) given in Table~\ref{tab:alt} and taken from~\cite{HTW}.

\begin{table}[htpb!]
\begin{center}
\begin{tabular}{|r|r|r|r|r|r|r|r|r|r|r|r|r|r|r|}
\hline    
crossings & 3 & 4 & 5 & 6 & 7 & 8 & 9 & 10 & 11 & 12 & 13 & 14 & 15 & 16  \\
\hline
\# alt. knots & 1 & 1 & 2 & 3 & 7 & 18 & 41 & 123 & 367 & \num{1288} & \num{4878} & \num{19536} &
\num{85263} & \num{379799} \\
\hline
\end{tabular}
\end{center}
\caption{Alternating knot counts, up to mirror image.}
\label{tab:alt}
\end{table}


After computing the $V_1$ and $V_2$ polynomials in the following range of knots,
we observed the following.

\begin{proposition}
\label{prop.pos}
For all alternating knots with $\leq 16$ crossings, we have
\be
\label{V1pos}
V_1(t,-q), V_2(t,-q) \in \BZ_{\geq 0}[t^{\pm 1},q^{\pm 1}] \,.
\ee
\end{proposition}

The above positivity fails for $V_3(t,-q)$ and $V_4(t,-q)$ already both for the
$3_1$ and the $4_1$ knots.

\begin{question}
\label{que.2}
Is this an accident of knots with low number of crossings or a hint of a
relation of $V_1$ and $V_2$ with some categorification theory?
\end{question}


\section{Computing the $V_n$-polynomials}
\label{sec.compute}

A priori, the polynomial invariant of long knots based on an $R$-matrix on a
$d$-dimensional vector space is a state sum of $c\cdot d^{2c-1}$ terms where $c$ is the
number of crossings of a planar projection of a knot, and in the case of
the $V_n$ polynomials, $d=4n$. Even though the summand is sparse, a direct
computation of the $V_2$ polynomial for knots with 8 crossings is unfeasible. 
A key observation is that every polynomial invariant of long knots of~\cite{GK:multi}
is given as a state sum and hence has a natural local tangle version.
The locality property of this polynomial is very important for its efficient
computation, an idea that is highlighted time and again in the work of Bar-Natan and
van der Veen (see e.g., ~\cite{BNV:API,BNV:PG}). In this section, we describe
the algorithm we used to compute the data presented in this paper and published
in \cite{GS:VnData}. A copy of our code for the computation of
$V_n$-polynomials ($n=1,2,3,4$), along with instructions, can also be found in
\cite{GS:VnData}. 

\subsection{Overview}
\label{subsec.algover}

Given a planar diagram of an oriented knot, the computation of the $V_n$-polynomial
is assembled from the following parts:
\begin{enumerate}[(I)]
\item
  Convert the planar diagram into a planar diagram of a corresponding oriented long
  knot with up-pointing crossings, define a height function on it and record the
  rotation number of each arc;
\item
  To each crossing and arc with nonzero rotation number, attach the $R$-matrix in
  consideration as a $4$-tensor, obtaining a tensor network on the shadow graph of
  the oriented long knot diagram;
\item
  Tensor contract the resulting tensor network, obtaining the $V_n$-polynomial.
\end{enumerate}
The real computation happens in Part~(III), which is exactly where the optimization
takes place. 

In terms of the above decomposition, Part~(I) is explained in
Section~\ref{subsub.longknot}, Part~(II) is explained by
Section~\ref{subsub.tensorcontract} and the first half of Section~\ref{subsec.compute},
and Part~(III) is explained by a mixture of Section~\ref{subsub.tensorcontract},
\ref{subsec.compute} and \ref{subsec.algorithm}.

\subsection{Terminology}
\label{subsec.term}

We introduce the terminology that will be used later in Section~\ref{subsec.compute}
and \ref{subsec.algorithm}. 

\subsubsection{Long knots and rotation numbers}
\label{subsub.longknot}

The $V_n$-polynomials, as defined in ~\cite{GK:multi}, are computed from the
\textit{oriented long knot diagrams} of knots. Given a planar diagram of an
oriented knot, one can obtain its corresponding oriented long knot diagram by the
following procedure:
\begin{itemize}
\item
  choose an arc in the planar diagram,
\item
  apply a stereographic projection so that the chosen arc bounds the face containing
  the infinity,
\item
  cut the chosen arc open, obtaining two open strands, and 
\item
  pull the two open strands to the infinity in two opposite directions without
  creating any self intersection.
\end{itemize}

As there are multiple choices of the arcs to cut, a knot may correspond to multiple
long knot diagrams. However, the $V_n$-polynomials are invariant under this choice
\cite{GHKST:link}, hence they are indeed invariants of oriented knots (as opposed
to oriented long knots).

In our settings, we pull the out-pointing strand vertically upward, and the other
strand to the opposite direction. By applying an ambient isotopy, we further require
that the segments of arcs always point upward in sufficiently small neighborhoods
of each crossing. The vertical axis now gives us a natural height function, which
restricts to a Morse function along the long knot diagram, whose local maxima
and minima occur only on the arcs and never at the crossings. We can now define
the rotation number associated to each arc in the oriented long knot diagram. 

\begin{definition}
The \textit{rotation number} of an arc $\mathcal{A}$ in an oriented long knot
diagram equipped with a height function $h$ as described above, is the integer
\begin{equation*}
  \sum_{p\in \{dh = 0 \} \cap \mathcal{A}} (-1)^{\delta_h(p)} \varepsilon_\mathcal{A} (p),
\end{equation*}
where $\delta_h(p) = 1$ if $p$ is a local maximum of $h$ and $0$ otherwise;
$\varepsilon_\mathcal{A} (p) = 1$ if the oriented arc $\mathcal{A}$ points to the
right at $p$, and $0$ otherwise. 
\end{definition}

\begin{figure}[htpb!]
\centering
\scalebox{.6}{$\vcenter{\hbox{\includegraphics{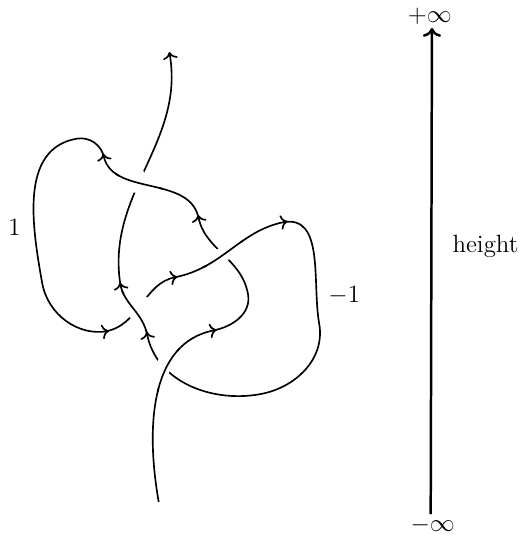}}}$
}
\caption{The long knot diagram corresponding to the $4_1$ knot, with nonzero
  rotation numbers labeled.}
\label{f.long4-1}
\end{figure}

See Figure~\ref{f.long4-1} for an example of oriented long knot diagram with
rotation numbers labeled. 

The definitions above are equivalent to those in \cite{BNV:API}. Note
that for most knots, different height functions (in the sense that the local maxima
and minima occur on different arcs) can be associated when we move their long knot
diagrams around with ambient isotopies. Hence the rotation numbers are not uniquely
determined by the oriented long knot diagram, but by the diagram together with the associated
height function. 

Given an oriented long knot diagram, one can always define a height function on it
so that the rotation number of each arc is either $\pm 1$ or $0$ by the following
procedure:
\begin{enumerate}[(i)]
  \item Start from the entrance strand.
  \item Walk along the arcs following the orientation while keep increasing
    height, stacking up-pointing crossings upon each other until having to turn
    around to connect the currently standing strand with an open strand attached
    to a lower crossing walked past previously. 
  \item Aftering connecting, continue to walk along the arc connected and repeat (ii). 
\end{enumerate}
The function \texttt{Rot[]} from \cite{BNV:API} computes the rotation numbers of
oriented long knot diagrams as if their height functions are defined as above. We
would like to remark that \texttt{Rot[]} does not work for links with multiple
components, and we provided a function \texttt{LinkRot[]} that is compatible with
links in the code shared in \cite{GS:VnData}.

\subsubsection{Tensor contraction}
\label{subsub.tensorcontract}

We now briefly review some standard terminology and results of tensor contraction.
For an extensive and graphical introduction to this topic, we refer the readers to
\cite{JC:tensornet}.

\begin{definition}
  Let $\mathcal{R}$ be a ring. An \textit{$n$-tensor} $T$ over $\mathcal{R}$ is a tuple
  $(T_{i_1,\dots,i_n})_{(i_1,\dots,i_n)\in \mathcal{S}}$ where $T_{i_1,\dots,i_n}\in \mathcal{R}$
  and the index set $\mathcal{S}$ is of the form 
\begin{equation*}
  \mathcal{S} = \left\{ 1,\dots,m_1 \right\}\times \dots\times \left\{
    1,\dots, m_n \right\} \subset \mathbb{N}^n,
\end{equation*}
We call the integer $m_k$ ($k\in \left\{ 1,\dots,n \right\}$) the dimension
of the \textit{leg} $k$ (an alias of the $k$-th index $i_k$) of the tensor $T$. 
\end{definition}

\begin{definition}
A \textit{tensor network} is a (multi)graph where each vertex is associated with
a tensor whose legs are bijectively associated with the ends of the edges incident
to the vertex. Each edge has two ends, which are associated with two different
legs of one or two tensors, and we further require the legs associated with the
ends of a same edge to have the same dimension.
\end{definition}

Given an $n$-tensor $T$ and an $n'$-tensor $T'$ over a common ring $\mathcal{R}$,
if the dimensions of leg $k$ of $T$ and leg $k'$ of $T'$ are equal, we can
\textit{contract} $T$ and $T'$ along the pair of
legs $(k,k')$ to obtain an $(n+n'-2)$-tensor $T''$, defined by
\begin{equation*}
T''_{i_1,\dots,\widehat{i_k},\dots,i_n,j_1,\dots,\widehat{j_{k'}},\dots,j_{n'}}
\coloneqq \sum_{\substack{i_k = j_{k'} \in \left\{ 1,\dots,m_k \right\}
}} T_{i_1,\dots,i_k,\dots,i_n}T'_{j_1,\dots,j_{k'},\dots,j_{n'}}
\end{equation*}
where the hats indicate that the corresponding indices are deleted. More generally,
if legs $k_1,\dots,k_s$ of $T$ have the same dimension respectively as legs
$k'_1,\dots,k'_s$ of $T'$, we can contract $T$ and $T'$ along the $s$ pairs of legs
$(k_1,k'_1),\dots, (k_s,k'_s)$ to obtain an $(n+n' - 2s)$-tensor $T''$ defined by
\begin{equation}
\label{tensorcontract}
T''_{i_1,\dots,\widehat{i_{k_1}},\dots,
  \widehat{i_{k_s}},\dots,i_n,j_1,\dots,\widehat{j_{k'_1}},\dots,
  \widehat{j_{k'_s}},\dots,j_{n'}}
\coloneqq \sum_{
\tiny
\substack{i_{k_1} = j_{k'_1} \in \{ 1,\dots,m_{k_1}\}
\\ \vdots\\
i_{k_s} = j_{k'_s} \in \{ 1,\dots,m_{k_s} \}
\normalsize
}} T_{i_1,\dots,i_n}T'_{j_1,\dots,j_{n'}}
\end{equation}
Similarly we can define contractions of more than two tensors at once.

It is clear from \eqref{tensorcontract} that, it demands $\mathcal{O}(c\cdot m_{k_1}\cdots m_{k_s})$ additions and multiplications to compute a single
entry in a tensor resulted from contracting $c$ tensors along $s$ pairs of legs
with common dimensions $m_{k_1},\dots,m_{k_s}$ at once. 

Let $m^{(l)}_1,\dots, m^{(l)}_{n^{(l)}}$ ($l\in \{1,\cdots,c\}$) denote the dimensions
of the legs of the $c$ tensors involved in the contraction, there are then
\begin{equation*}
  \frac{\prod_{l=1}^c \left(m^{(l)}_1\cdots m^{(l)}_{n^{(l)}}\right)}{\left(m_{k_1}
      \cdots m_{k_s}\right)^2}
\end{equation*}
entries in the resulting tensor. Therefore the total number of additions and multiplications required for
the tensor contraction described above is in 
\begin{equation*}
 \mathcal{O}\left( c\cdot (m_{k_1}\cdots m_{k_s})\right)\cdot \frac{\prod_{l=1}^c \left(m^{(l)}_1\cdots
      m^{(l)}_{n^{(l)}}\right)}{\left(m_{k_1}\cdots m_{k_s}\right)^2}
  = \mathcal{O}\left(c\cdot \frac{\prod_{l=1}^c \left(m^{(l)}_1\cdots m^{(l)}_{n^{(l)}}\right)}{m_{k_1}
    \cdots m_{k_s}}\right).
\end{equation*}
In other words, the number of additions and multiplications for a single tensor contraction can be
quantified by
\begin{equation}
\label{TCC}
\text{Number of tensors involved}\cdot\frac{\text{Product of dimensions of all legs of tensors
    involved}}{\text{Product of dimensions of all pairs of legs contracted}}.
\end{equation}
Note that the dimensions in the denominator of \eqref{TCC} count only once for
each pair (which involves two legs) of legs contracted. Note also that the computational
complexity of doing the additions and multiplications in ring 
$\mathcal{R}$ is not quantified by \eqref{TCC}. In the computation of $V_n$-polynomials,
$\mathcal{R}$ will be the ring $\mathbb{Z}[t^{\pm 1}, q^{\pm 1}]$.

\begin{example}
A matrix can be seen as a $2$-tensor, with its number of rows and number of
columns being the dimensions of its two legs. Let $ A$ be an $n\times m$-matrix
and $ B$ be an $m\times k$-matrix, computing the matrix multiplication
$ A B$ is the same as computing the contraction of the corresponding two
$2$-tensors along the pair of legs of dimension $m$. According to \eqref{TCC},
the time complexity of this computation is thus
\begin{equation*}
2\cdot \frac{n\cdot m \cdot m\cdot k }{m} = 2nmk .
\end{equation*}
This recovers the naive time complexity of matrix multiplications. 
\end{example}

\begin{example}
\label{TCC-d}
Let all legs be of dimension $d$ in this example. The number of additions and multiplications for contracting
$c$ tensors at once along $s$ pairs of legs, where the tensors have
$n_1,\dots, n_c$ legs respectively, is
\begin{equation*}
  c\cdot \frac{d^{n_1}\cdots d^{n_c}}{d^s} = c\cdot d^{n_1 + \cdots + n_c - s}.
\end{equation*}
\end{example}

\subsection{Process of computation}
\label{subsec.compute}

Before delving into the discussion of the algorithms, we give a general description
of how the $V_n$-polynomials are computed without touching the theoretical parts,
for which we kindly refer the readers to \cite{GK:multi,Ka:longknots}. 

Let $R\in \Aut(V\otimes V)$ be the $R$-matrix of the $V_n$-polynomial we wish to
compute, then $\dim V = 4n$ and we fix a basis $\mathcal{B} \coloneqq
\{e_1,\dots,e_{4n}\}$ of $V$. Let $(R^{\pm 1})_{e_i\otimes e_j} ^{e_k\otimes e_l}$ be
the matrix entry of $R^{\pm 1}$ corresponding to the basis element $e_i\otimes e_j$
in the domain and $e_k\otimes e_l$ in the codomain, thus
$\{(R^{\pm 1})_{e_i\otimes e_j} ^{e_k\otimes e_l}\}_{(i,j,k,l)}$ becomes a $4$-tensor
with dimension $4n$ for all four legs. 

Now we start from an oriented long knot diagram with a height function equipped,
as is set up in Section~\ref{subsub.longknot}, and label all arcs with $a_i$'s.
We locally replace the crossings with the $4$-tensors as the following: 
\begin{equation*}
  \vcenter{\hbox{\includegraphics{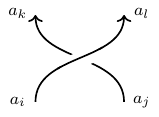}}} 
    \quad \mapsto \quad 
    \vcenter{\hbox{\includegraphics{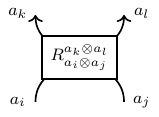}}},
    \qquad
    \vcenter{\hbox{\includegraphics{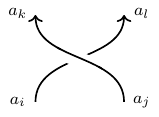}}}
    \quad \mapsto \quad 
    \vcenter{\hbox{\includegraphics{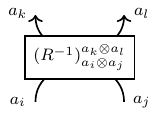}}}
\end{equation*}
Note that there are only these two kinds of crossings due to our requirement that
all crossings point up. Additionally, we put a $2$-tensor
$((\widetilde{R})^{-1} \eta)_{a_i}^{a_i}$ or $(\varepsilon(
\widetilde{R^{-1}})^{-1})_{a_i}^{a_i}$ at arcs $a_i$ with rotation number $1$ or
$-1$ respectively; for conciseness we refer the readers again to \cite{GK:multi}
for the definitions of $(\widetilde{R})^{-1} \eta$ and $\varepsilon(
\widetilde{R^{-1}})^{-1}$, and, for the sake of complexity analysis, mention only
the fact that under our choice of the basis of $V$ they will be diagonal matrices
with only $\pm1$ along the diagonals. 

The above procedure gives us a tensor network with two free ends. Contracting it as
all $a_i$'s run through the basis of $V$ except for the entrance and exit strands,
we obtain a $2$-tensor whose corresponding matrix is a scalar
multiplication \cite{GHKST:link}. The eigenvalue of this
matrix is the desired $V_n$-polynomial. 

\begin{figure}[htpb!]
\centering
    \begin{equation*}
\vcenter{\hbox{\includegraphics{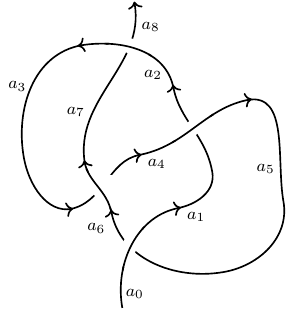}}}
    \quad\rightsquigarrow \quad 
    \vcenter{\hbox{\includegraphics{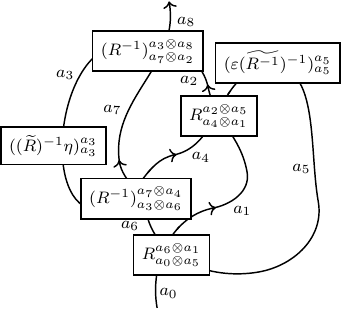}}}
\end{equation*}
\caption{The tensor network resulted from an oriented long knot diagram
      of the $4_1$ knot.}
\label{f.tensornet4-1}
\end{figure}

\begin{example}
  \label{ex.4-1-v}
  Figure~\ref{f.tensornet4-1} illustrates the procedure of obtaining the
  corresponding tensor network on an oriented long knot diagram of the $4_1$ knot.
  With the symbols in Figure~\ref{f.tensornet4-1}, the $V_n$-polynomial of the $4_1$
  knot can be expressed as the following sum:
  \begin{equation*}
    \sum_{\substack{a_1,\dots,a_7\in \mathcal{B}\\ a_0 = a_8 = e_1}}
    R_{a_0\otimes a_5}^{a_6\otimes a_1} \cdot 
    ((\widetilde{R})^{-1} \eta)_{a_3}^{a_3}\cdot  
    (R^{-1})_{a_3\otimes a_6}^{a_7\otimes a_4} \cdot 
    R_{a_4\otimes a_1}^{a_2\otimes a_5} \cdot
    (\varepsilon( \widetilde{R^{-1}})^{-1})_{a_5}^{a_5} \cdot
    (R^{-1})_{a_7\otimes a_2}^{a_3\otimes a_8}.
  \end{equation*}
  Recall that $\mathcal{B} = \{e_1,\dots,e_{4n}\}$ is the basis of $V$.
\end{example}

Since the formula, as in Example~\ref{ex.4-1-v}, is a finite sum, we can compute
it piece by piece in an arbitrary order; the time complexity, however, is dependent
on that order. This is where the optimization algorithm comes into play. 

\subsection{The algorithm}
\label{subsec.algorithm}

To optimize the computation, we need to find a sequence of tensor contractions
such that the time complexity of performing the contractions following said sequence
is as small as possible. For general tensor networks, the optimization problem of
finding such a sequence is well-known to be NP-complete \cite{LSR:NPtensor}. However,
the tensor networks here in our discussion, arising from knots, are much less
complicated, hence we are able to reduce the complexity to a satisfactory level
without tackling the NP-hardness. 

Before we start, we set the labels of the entrance and exit strands to $e_1$ to
avoid repetitive computation. We also ignore the $2$-tensors we attached to arcs
with nonzero rotation numbers in the following analysis, since they only represent
changes of signs, which are trivial to compute. 

We now have a tensor network with only $4$-tensors and two $3$-tensors, and all
legs with the same dimension $4n$. Applying $\log_{4n}$ to \eqref{TCC} gives the
following description of the number $N$ of additions and multiplications for each contraction:
\begin{equation}
\label{LogTCC}
\begin{aligned}
\log_{4n}N =&\ \log_{4n}(\text{number of tensors involved})\\ 
& +  \#\{\text{legs of tensors involved}\} - \#\{\text{pairs of legs contracted}\},
\end{aligned}
\end{equation}
where $\# S$ stands for the cardinality of a set $S$. More concisely, using the
same convention as in Example~\ref{TCC-d} with $d = 4n$, we can write \eqref{LogTCC} as 
\begin{equation*}
  \log_{d}N = \log_d c + n_1 + \cdots + n_c - s.
\end{equation*}

We call $ n_1 + \cdots + n_c - s $ the \textit{local
  contraction width} of a contraction operation. Given a sequence of contractions,
its \textit{contraction width} $w$ is the maximum of the local contraction widths
of all contractions within that sequence, and its \textit{multiplicity} $m$ is the
number of contractions whose local contraction widths are equal to $w$. 
The total number of additions and multiplications for computing the $V_n$-polynomial is thus
  $\mathcal{O}(D\cdot m\cdot d^{w})$, where $D$ is the number of crossings in the knot diagram (hence $c\leq D$). With the knot diagram fixed, we may simply use the pair $(w,m)$,
under the lexicographical ordering, to estimate the overall difficulty of computation. 

\begin{figure}[htpb!]
  \centering
  $\vcenter{\hbox{\includegraphics{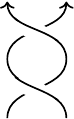}}}$
\caption{A bigon in a knot diagram.}
\label{f.bigon}
\end{figure}

Therefore, our task is to find a sequence whose pair $(w,m)$ is satisfactorily small
for us to perform the computation. A previously known approach is to prioritize
contracting bigons, i.e. two tensors sharing two legs (see Figure~\ref{f.bigon}),
as is described in \cite[Section 5.5]{MST:QES}. In tensor networks obtained from
knot diagrams, contracting bigon is the simplest move one can take, whose local
contraction width is generically $6$, and $5$ if one of the crossing contains the
entrance or exit strand. However, this approach lacks means of optimization in the
middle of the procedure, and especially struggles when dealing with knots with few
bigons in its diagram, an extreme example being the $9_{40}$ knot whose initial
diagram does not have any bigon. Hence we take a different approach, which can be
best summarized as a local minimum search: 

\begin{description}
  \item[Input] An oriented knot diagram.
  \item[Procedure] 
\begin{enumerate}[(A)]
\item Generate a list $\mathcal{L}$ of oriented long knot diagrams obtained by
  cutting each arc in the knot diagram open.
\item For each diagram in $\mathcal{L}$, remove the entrance and exit strands,
  consider its corresponding tensor network, obtain a sequence of contraction
  operations by the following and record its contraction width $w$ and multiplicity
  $m$:
\begin{description}
  \item[While] There are edges in the tensor network remaining,
  \begin{enumerate}
  \item
    Find all possible contraction operations of two tensors with common edges;
  \item
    Evaluate the local contraction width for each operation;
  \item
    Update the tensor network as if the contraction with the minimum local
    contraction width is performed (if there are multiple such contractions,
    choose the first one in the list).
  \end{enumerate}
\end{description}
\end{enumerate}
\item[Output] The diagram in $\mathcal{L}$ with the minimal $(w,m)$ under the
  lexicographical ordering.
\end{description}

With the output diagram, we reproduce the sequence obtained in Step~(B) and
perform the actual computation following it.

We implemented the above algorithm in \texttt{Mathematica} to perform the computation, and our code is available at \cite{GS:VnData}. It took about a month for our algorithm to compute the $V_2$-polynomial for all $\num{1388705}$ knots of $16$ crossings in 300 threads.

For the extreme example, the $9_{40}$ knot, the algorithm of collapsing bigons
took about $3$ minutes to compute its $V_2$-polynomial, while our algorithm above
took $30$ seconds. More generally, for the $40$ randomly chosen knots of $15$
crossings in Table~\ref{tab:randomknots}, the algorithm of collapsing bigons took
about $105$ minutes to compute their $V_2$-polynomials, while our algorithm took
$25$ minutes.

\begin{table}[htpb!]
\tiny
\begin{equation*}
  \begin{aligned}
    15n3881, && 15n99297,&& 15a44733,&& 15n11237, &&15n70455, &&15n83510,&&15a36536,
    && 15n157698, \\
    15a83512, &&15a38060,&& 15n100625, &&15n45279, &&15a73940, &&15n120336,&& 15n154365,
    &&15a38391, \\
15a20828,&& 15a55989, &&15a42521,&& 15n126571, &&
15a51133, &&15n5288,&& 15a58690, &&15a28867,\\
15a82438,&& 15a63595, &&15n89407, &&15n136788, &&15a84489, &&15n160922, &&15a11824,
&&15n165397,\\
15a12009,&& 15a26849,&& 15a60295, &&
15a36928, &&15n16661,&& 15n27741, &&15n67039, &&15n151938.
  \end{aligned}
\end{equation*}
\normalsize
\normalsize
\caption{The $40$ randomly chosen knots of $15$ crossings for performance
  comparison.}
\label{tab:randomknots}
\end{table}


\subsection*{Acknowledgements}

The authors wish to thank Dror Bar-Natan, Nathan Dunfield, Rinat Kashaev,
Ben-Michael Kohli, Ciprian Manolescu, Mingde Ren, Roland van der Veen and
Cai Zeng for useful conversations. This material is based upon work supported by the National Science Foundation under Grant No. DMS-2424139 while the second author was in residence at the Simons Laufer Mathematical Sciences Institute in Berkeley,
California, during the Spring 2026 semester. The second author was also supported by the National Science Foundation under Grant No. DMS-2303572 during the same semester.

\subsection*{Conflict of interest}
To the best of our knowledge, the authors are not aware of any conflict
of interest.


\clearpage
\appendix

\section{$V_2$-equivalence classes of knots with at most 14 crossings}
\label{app.pairs}

In this section we list the $V_2$-equivalence classes (up to mirror image)
of knots with at most 14 crossings. As perhaps expected, the equivalence classes
involve knots with the same number of crossings. Overline means mirror image. 
The counts are given in Table~\ref{tab:3}.



Below, we indicate the flavor of each equivalence class (defined in
Equation~\eqref{3flavors}) by the corresponding number in the superscript.
There are 387 tuples in total, 293 being tight + thin, 57 tight + thick,
and 37 loose + thick. 
  
First, we give the 3 pairs with 12 crossings.

\begin{tiny}
\be
\label{pairs12}
\begin{aligned}
(12n364, \,\,\, \overline{12n365})^1  &&
(12n421, \,\,\, \overline{12n422})^1  &&
(12n553, \,\,\, \overline{12n556})^1  
\end{aligned}
\ee
\end{tiny}

Next, we give the 50 pairs of knots with 13 crossings.

\begin{tiny}
\be
\label{pairs13}
\begin{aligned}
(13a141, \,13a142)^1 &&
(13a150, \,13a154)^1 &&
(13a199, \,13a204)^1 &&
(13a316, \,13a349)^1 \\
(13a906, \,13a916)^1 &&
(13a967, \,13a1059)^1 &&
(13a1114, \,13a1143)^1 &&
(13a1123, \,13a1394)^1 \\
(13a1126, \,13a1163)^1 &&
(13a1781, \,13a1816)^1 &&
(13a1813, \,13a1831)^1 &&
(13a1990, \,13a2020)^1 \\
(13a1991, \,13a2021)^1 &&
(13a1993, \,13a2024)^1 &&
(13a1995, \,13a2006)^1 &&
(13a1996, \,13a2022)^1 \\
(13a2720, \,13a2727)^1 &&
(13a2800, \,13a2801)^1 &&
(13a2802, \,13a2808)^1 &&
(13a3238, \,13a3254)^1 \\
(13n370, \,13n373)^1 &&
(13n371, \,13n374)^1 &&
(13n372, \,13n375)^3 &&
(13n403, \,13n415)^1 \\
(13n404, \,13n416)^1 &&
(13n405, \,13n417)^2 &&
(13n406, \,13n418)^1 &&
(13n407, \,13n419)^2 \\
(13n534, \,13n549)^1 &&
(13n535, \,13n550)^1 &&
(13n536, \,13n551)^3 &&
(13n874, \,13n949)^1 \\
(13n875, \,13n950)^1 &&
(13n876, \,13n951)^3 &&
(13n1129, \,13n1130)^1 &&
(13n1359, \,13n1360)^1 \\
(13n1653, \,13n1683)^3 &&
(13n1654, \,13n1684)^1 &&
(13n1655, \,13n1685)^2 &&
(13n1893, \,13n2098)^1 \\
(13n1894, \,13n2099)^1 &&
(13n2184, \,13n2228)^1 &&
(13n2185, \,13n2229)^1 &&
(13n2186, \,13n2201)^3 \\
(13n2205, \,13n2250)^3 &&
(13n2930, \,13n2977)^3 &&
(13n2933, \,\overline{13n2956})^1 &&
(13n2934, \,\overline{13n2954})^3 \\
(13n2937, \,13n2955)^1 &&
(13n3510, \,13n3517)^3 &&
\end{aligned}
\ee
\end{tiny}

Next, we give the 333 pairs and one triple of knots with 14 crossings. 

\begin{tiny}
\begin{subequations}
\be
\label{pairs14a}
\begin{aligned}
(14a34, \,14a35)^1 &&
(14a43, \,14a48)^1 &&
(14a96, \,14a103)^1 &&
(14a195, \,14a228)^1 \\
(14a454, \,14a458)^1 &&
(14a516, \,14a520)^1 &&
(14a518, \,14a592)^1 &&
(14a519, \,14a593)^1 \\
(14a533, \,14a550)^1 &&
(14a534, \,14a549)^1 &&
(14a608, \,14a617)^1 &&
(14a631, \,14a635)^1 \\
(14a675, \,14a734)^1 &&
(14a717, \,14a735)^1 &&
(14a718, \,14a736)^1 &&
(14a780, \,14a786)^1 \\
(14a989, \,14a1017)^1 &&
(14a1044, \,14a1166)^1 &&
(14a1047, \,14a1170)^1 &&
(14a1048, \,14a1168)^1 \\
(14a1268, \,14a1362)^1 &&
(14a1436, \,14a1437)^1 &&
(14a1445, \,14a1449)^1 &&
(14a1488, \,14a1493)^1 \\
(14a1522, \,14a1532)^1 &&
(14a1660, \,14a1724)^1 &&
(14a1767, \,14a1860)^1 &&
(14a2083, \,14a2113)^1 \\
(14a2205, \,14a2215)^1 &&
(14a2245, \,14a2254)^1 &&
(14a2253, \,14a2256)^1 &&
(14a2567, \,14a2573)^1 \\
(14a2609, \,14a2618)^1 &&
(14a2617, \,14a2620)^1 &&
(14a3400, \,14a3433)^1 &&
(14a3402, \,14a3434)^1 \\
(14a3403, \,14a3432)^1 &&
(14a3405, \,14a3435)^1 &&
(14a3407, \,14a3436)^1 &&
(14a3408, \,14a3437)^1 \\
(14a3409, \,14a3438)^1 &&
(14a3412, \,14a3416)^1 &&
(14a3419, \,14a3439)^1 &&
(14a3539, \,14a4433)^1 \\
(14a4041, \,14a4998)^1 &&
(14a4046, \,14a4140)^1 &&
(14a4147, \,14a4939)^1 &&
(14a4264, \,14a4871)^1 \\
(14a4901, \,14a5698)^1 &&
(14a4904, \,14a5077)^1 &&
(14a6467, \,\overline{14a6614})^1 &&
(14a6468, \,14a6620)^1 \\
(14a7190, \,14a7247)^1 &&
(14a7193, \,14a7216)^1 &&
(14a7194, \,14a7218)^1 &&
(14a7196, \,14a7269)^1 \\
(14a7199, \,14a7202)^1 &&
(14a7200, \,14a7249)^1 &&
(14a7205, \,14a7271)^1 &&
(14a7207, \,14a7272)^1 
\end{aligned}
\ee

\be
\label{pairs14b}
\begin{aligned}
(14a7209, \,14a7276)^1 &&
(14a7210, \,14a7275)^1 &&
(14a7213, \,14a7279)^1 &&
(14a7215, \,14a7698)^1 \\
(14a7217, \,14a7270)^1 &&
(14a7219, \,14a7260)^1 &&
(14a7220, \,14a7277)^1 &&
(14a7258, \,14a7263)^1 \\
(14a7262, \,14a7273)^1 &&
(14a7264, \,14a7274)^1 &&
(14a7399, \,14a7449)^1 &&
(14a7446, \,14a7477)^1 \\
(14a7478, \,14a8066)^1 &&
(14a7527, \,14a7598)^1 &&
(14a8016, \,14a8106)^1 &&
(14a8017, \,14a8107)^1 \\
(14a8023, \,14a8321)^1 &&
(14a8096, \,14a8115)^1 &&
(14a8829, \,14a9027)^1 &&
(14a9707, \,14a9711)^1 \\
(14a10115, \,14a10156)^1 &&
(14a10116, \,\overline{14a10142})^1 &&
(14a10120, \,14a10170)^1 &&
(14a10163, \,14a10185)^1 \\
(14a10405, \,14a10410)^1 &&
(14a10406, \,14a10455)^1 &&
(14a10407, \,14a10434)^1 &&
(14a10408, \,14a10412)^1 \\
(14a10411, \,14a10453)^1 &&
(14a10413, \,14a10440)^1 &&
(14a10414, \,14a10417)^1 &&
(14a10415, \,14a10435)^1 \\
(14a10437, \,14a10452)^1 &&
(14a10439, \,14a10456)^1 &&
(14a10443, \,14a10451)^1 &&
(14a10853, \,14a11544)^1 \\
(14a11322, \,14a11553)^1 &&
(14a11793, \,14a12325)^1 &&
(14a12333, \,14a12334)^1 &&
(14a12335, \,14a12344)^1 \\
(14a12336, \,14a12337)^1 &&
(14a12816, \,14a12833)^1 &&
(14a12817, \,14a12834)^1 &&
(14a12868, \,14a12876)^1 \\
(14a12869, \,14a12877)^1 &&
(14a13431, \,14a13433)^1 &&
(14a13434, \,14a13436)^1 &&
(14a13473, \,14a13475)^1 \\
(14a13476, \,14a13478)^1 &&
(14n179, \,14n182)^1 &&
(14n180, \,14n183)^1 &&
(14n181, \,14n184)^3 \\
(14n212, \,14n225)^1 &&
(14n213, \,14n226)^1 &&
(14n214, \,14n227)^2 &&
(14n215, \,14n228)^1 \\
(14n216, \,14n229)^3 &&
(14n364, \,14n386)^1 &&
(14n365, \,14n387)^1 &&
(14n366, \,14n388)^3 \\
(14n732, \,14n809)^1 &&
(14n733, \,14n810)^1 &&
(14n734, \,14n811)^3 &&
(14n1366, \,14n1393)^1 \\
(14n1369, \,14n1394)^3 &&
(14n1370, \,14n1395)^2 &&
(14n1373, \,14n1396)^1 &&
(14n1374, \,14n1397)^1 \\
(14n1377, \,14n1398)^1 &&
(14n1378, \,14n1399)^1 &&
(14n1379, \,14n1400)^1 &&
(14n1380, \,14n1401)^3 \\
(14n1691, \,14n1703)^1 &&
(14n1692, \,14n1704)^1 &&
(14n1693, \,14n1705)^3 &&
(14n1697, \,14n1945)^1 \\
(14n1698, \,14n1946)^1 &&
(14n1699, \,14n1947)^2 &&
(14n1700, \,14n1948)^1 &&
(14n1701, \,14n1949)^1 \\
(14n1702, \,14n1950)^2 &&
(14n1752, \,14n1753)^1 &&
(14n1760, \,14n1761)^1 &&
(14n1762, \,14n1839)^2 \\
(14n1763, \,14n1840)^1 &&
(14n1764, \,14n1841)^2 &&
(14n1765, \,14n1842)^1 &&
(14n1766, \,14n1843)^2 \\
(14n1767, \,14n1774)^1 &&
(14n1768, \,14n1775)^1 &&
(14n1769, \,14n1834)^1 &&
(14n1770, \,14n1835)^1 \\
(14n1771, \,14n1836)^2 &&
(14n1772, \,14n1837)^1 &&
(14n1773, \,14n1838)^2 &&
(14n2007, \,14n2032)^2 \\
(14n2039, \,14n2042)^3 &&
(14n2148, \,14n2372)^1 &&
(14n2149, \,14n2373)^1 &&
(14n2150, \,14n2374)^3 \\
(14n2294, \,14n2375)^1 &&
(14n2295, \,14n2376)^1 &&
(14n2296, \,14n2377)^2 &&
(14n2297, \,14n2378)^1 \\
(14n2298, \,14n2379)^1 &&
(14n2299, \,14n2380)^2 &&
(14n2420, \,14n4659)^3 &&
(14n2423, \,14n5868)^2 \\
(14n3406, \,14n3407)^1 &&
(14n3418, \,14n3823)^1 &&
(14n3421, \,14n3824)^1 &&
(14n3422, \,14n3825)^3 \\
(14n3443, \,14n3835)^1 &&
(14n3444, \,14n3836)^1 &&
(14n3445, \,14n3837)^2 &&
(14n3448, \,14n3829)^1 \\
(14n3449, \,14n3830)^1 &&
(14n3450, \,14n3831)^2 &&
(14n4561, \,14n4562)^3 &&
(14n4577, \,14n4583)^2 \\
(14n4578, \,14n4584)^1 &&
(14n4579, \,14n4585)^2 &&
(14n4619, \,14n4625)^3 &&
(14n4657, \,14n4665)^1 \\
(14n4658, \,14n4664)^1 &&
(14n4925, \,14n5085)^1 &&
(14n4926, \,14n4927)^1 &&
(14n4930, \,14n4931)^1 \\
(14n4932, \,14n5086)^1 &&
(14n4933, \,14n5087)^2 &&
(14n4934, \,14n4935)^2 &&
(14n4938, \,14n4939)^1 \\
(14n5756, \,14n5780)^3 &&
(14n5822, \,14n5852)^3 &&
(14n5828, \,14n6370)^2 &&
(14n5854, \,14n5862)^1 \\
(14n5855, \,14n5861)^1 &&
(14n7506, \,14n7559)^1 &&
(14n7565, \,14n7674)^3 &&
(14n7566, \,14n7673)^1 \\
(14n7567, \,14n7672)^1 &&
(14n7575, \,14n7675)^2 &&
(14n7576, \,14n7676)^1 &&
(14n7577, \,14n7677)^2 \\
(14n7580, \,14n7671)^3 &&
(14n7586, \,14n7678)^1 &&
(14n7592, \,14n7679)^2 &&
(14n7593, \,14n7680)^1 \\
(14n7594, \,14n7681)^2 &&
(14n7597, \,14n7682)^2 &&
(14n7598, \,14n7683)^1 &&
(14n7599, \,14n7684)^2 \\
(14n7602, \,14n7685)^2 &&
(14n7603, \,14n7686)^1 &&
(14n7604, \,14n7687)^2 &&
(14n7617, \,14n7628)^2 \\
(14n7618, \,14n7629)^1 &&
(14n7636, \,14n7688)^1 &&
(14n7637, \,14n7689)^2 &&
(14n7638, \,14n7690)^1 \\
(14n7639, \,14n7691)^2 &&
(14n8225, \,14n10806)^2 &&
(14n8226, \,14n10807)^1 &&
(14n8291, \,\overline{14n8293})^1 \\
(14n8648, \,\overline{14n8649})^1 &&
(14n8650, \,\overline{14n8651})^1 &&
(14n8696, \,\overline{14n8697})^1 &&
(14n9075, \,\overline{14n9076})^1 \\
(14n9139, \,\overline{14n9140})^1 &&
(14n9142, \,\overline{14n9143})^1 &&
(14n9395, \,\overline{14n9396})^1 &&
(14n9398, \,\overline{14n9399})^1 \\
(14n9455, \,\overline{14n9456})^1 &&
(14n9458, \,\overline{14n9459})^1 &&
(14n9686, \,\overline{14n9687})^1 &&
(14n10002, \,14n11740)^3 \\
(14n10502, \,14n11639)^1 &&
(14n10503, \,14n11641)^3 &&
(14n10504, \,14n11640)^1 &&
(14n11679, \,14n11981)^3 
\end{aligned}
\ee

\be
\begin{aligned}
\label{pairs14c}
(14n14119, \,14n14286)^2 &&
(14n14122, \,14n14288)^1 &&
(14n14123, \,14n14287)^2 &&
(14n14130, \,\overline{14n14216})^1 \\
(14n14131, \,\overline{14n14214})^2 &&
(14n14134, \,14n14215)^2 &&
(14n14135, \,\overline{14n14221})^1 &&
(14n14136, \,\overline{14n14219})^3 \\
(14n14139, \,14n14220)^1 &&
(14n14148, \,14n14150)^1 &&
(14n14149, \,14n14151)^1 &&
(14n14153, \,14n14327)^2 \\
(14n14154, \,14n14328)^1 &&
(14n14155, \,14n14329)^2 &&
(14n14156, \,\overline{14n14157})^1 &&
(14n14158, \,\overline{14n14159})^1 \\
(14n14162, \,14n14169)^1 &&
(14n14165, \,14n14292)^1 &&
(14n14177, \,14n14330)^3 &&
(14n14192, \,14n14334)^1 \\
(14n14196, \,14n14333)^1 &&
(14n14203, \,14n14340)^1 &&
(14n14204, \,14n14341)^1 &&
(14n14205, \,14n14342)^2 \\
(14n14208, \,14n15068)^2 &&
(14n14209, \,14n14224)^2 &&
(14n14210, \,14n14225)^1 &&
(14n14211, \,\overline{14n14223})^2 \\
(14n14226, \,14n14315)^1 &&
(14n14227, \,14n14335)^1 &&
(14n14228, \,14n14336)^1 &&
(14n14313, \,14n14319)^3 \\
(14n14318, \,14n14331)^3 &&
(14n14322, \,14n14332)^2 &&
(14n14502, \,\overline{14n14508})^1 &&
(14n14504, \,\overline{14n14506})^1 \\
(14n14509, \,\overline{14n14516})^1 &&
(14n14511, \,\overline{14n14513})^1 &&
(14n14589, \,14n14662)^1 &&
(14n14590, \,14n14663)^1 \\
(14n14654, \,14n14684)^1 &&
(14n14655, \,14n14685)^1 &&
(14n14656, \,14n14686)^2 &&
(14n14687, \,14n15694)^1 \\
(14n14688, \,14n15695)^1 &&
(14n14780, \,14n14893)^2 &&
(14n14786, \,14n14894)^1 &&
(14n14787, \,14n14895)^1 \\
(14n14788, \,14n14896)^3 &&
(14n14793, \,14n14897)^1 &&
(14n14804, \,\overline{14n14808})^1 &&
(14n14923, \,\overline{14n14924})^1 \\
(14n14925, \,\overline{14n14926})^1 &&
(14n14930, \,\overline{14n14931})^1 &&
(14n15022, \,\overline{14n15024})^1 &&
(14n15058, \,\overline{14n15059})^1 \\
(14n15062, \,\overline{14n15063})^1 &&
(14n15065, \,\overline{14n15066})^1 &&
(14n15083, \,\overline{14n15084})^1 &&
(14n15102, \,\overline{14n15103})^1 \\
(14n15105, \,\overline{14n15106})^1 &&
(14n15172, \,\overline{14n15173})^1 &&
(14n15179, \,\overline{14n15180})^1 &&
(14n15201, \,\overline{14n15202})^1 \\
(14n15204, \,\overline{14n15205})^1 &&
(14n15207, \,\overline{14n15208})^1 &&
(14n15227, \,\overline{14n15228})^1 &&
(14n15231, \,\overline{14n15232})^1 \\
(14n15257, \,\overline{14n15258})^1 &&
(14n15630, \,14n15975)^3 &&
(14n15727, \,14n15756)^1 &&
(14n15728, \,14n15757)^1 \\
(14n15729, \,14n15758)^2 &&
(14n16547, \,14n16669)^3 &&
(14n17934, \,\overline{14n17940})^1 &&
(14n17936, \,\overline{14n17938})^2 \\
(14n17939, \,14n17945)^1 &&
(14n17941, \,14n18015)^2 &&
(14n17942, \,14n17946)^1 &&
(14n17948, \,\overline{14n17962})^1 \\
(14n17949, \,\overline{14n17960})^2 &&
(14n17952, \,14n17961)^1 &&
(14n17986, \,14n18013)^1 &&
(14n17994, \,14n18016)^1 \\
(14n17997, \,14n18017)^1 &&
(14n17998, \,14n18018)^2 &&
(14n18005, \,14n18012)^2 &&
(14n18144, \,\overline{14n18146})^1 \\
(14n18207, \,\overline{14n18208})^1 &&
(14n18947, \,14n19733)^3 &&
(14n19744, \,14n19758)^2 &&
(14n20141, \,14n20174)^2 \\
(14n20142, \,14n20175)^2 &&
\end{aligned}
\ee

\be
\label{pairs14d}
\begin{aligned}
  (14n14212, \,14n14213, \,\overline{14n14222})^1
\end{aligned}
\ee
\end{subequations}
\end{tiny}




The knots in Tables~\eqref{pairs12}, ~\eqref{pairs13},
~\eqref{pairs14a}--\eqref{pairs14d} are Conway mutant, have equal HFK and
Khovanov Homology and Equation~\eqref{degVn} is an equality for $n=2$. The knots in
Tables~\eqref{pairs12}, ~\eqref{pairs13} have equal $V_1$, $V_2$, $V_3$ and $V_4$
polynomials and the ones in Tables ~\eqref{pairs14a}--\eqref{pairs14d} have equal
$V_1$, $V_2$ and $V_3$ polynomials with the exceptions in Table~\ref{tab:exceptional}.

Finally a comment about quasi-alternating knots. By a computer search, we 
confirmed that, out of the $329$ non-alternating
tight + thin knots listed above, all but the following $20$ are quasi-alternating: 
\begin{tiny}
\begin{equation}
  \label{nonQA}
  \begin{aligned}
  13n403 &&
  13n2098 &&
  14n212 &&
  14n1775 &&
  14n2378 &&
  14n3448 &&
  14n4925 &&
  14n5085 &&
  14n5854 && 
  14n5862 \\ 
  14n7506 && 
  14n7559 && 
  14n14135 && 
  14n14221 && 
  14n14149 && 
  14n14151 &&
  14n14162 && 
  14n14169 &&
  14n14165 &&
  14n14292 
  \end{aligned}
\end{equation}
\end{tiny}
Note that $14n14135$ is the only $\overline{\text{Kh}'}$-thick knot in \eqref{nonQA},
hence the rest are all candidates of homologically thin non-QA knots.

As a byproduct of our computer search, we confirmed that the knots $12n139$ and
$12n331$, previously considered two candidates of $12$-crossing homologically thin
non-QA knots in \cite{Jablan}, are in fact quasi-alternating, with the following
initial planar diagrams and crossings (bolded) to perform the first smoothing.

\begin{tiny}
\begin{equation}
\begin{aligned}
  12n139\colon & [(27, 16, 28, 17), (17, 26, 18, 27), \bm{(18, 10, 19, 9)},
  (1, 22, 2, 23), (12, 20, 13, 19), (8, 16, 9, 15), (6, 23, 7, 24), (4, 14, 5, 13),\\
  & (21, 2, 22, 3), (29, 11, 30, 10), (31, 21, 32, 20), (3, 1, 4, 32),
  (11, 31, 12, 30), (24, 7, 25, 8), (25, 29, 26, 28), (14, 6, 15, 5)]\\
  12n331 \colon & [(17, 11, 18, 10), (23, 13, 24, 12), (27, 18, 28, 19),
  (24, 2, 25, 1), (3, 23, 4, 22), (6, 20, 7, 19), (11, 4, 12, 5), \\
 & \bm{(13, 3, 14, 2)},
 (9, 17, 10, 16), (7, 26, 8, 27), (20, 25, 21, 26), (15, 9, 16, 8),
 (14, 22, 15, 21), (28, 5, 1, 6)]
  \end{aligned}
\end{equation}
\end{tiny}

We plan to share our code of certifying quasi-alternating knots in a future
publication, where we give an extended table of quasi-alternating links compared
to that in \cite{Jablan}.

  

\section{Values for Whitehead doubles and (2,1)-parallels}

The values of $V_{\mathrm{Wh}(K),2}(t,q)$ for the first three nontrivial knots
is given as follows, where $u=t+t^{-1}-q-q^{-1}$ is as in~\eqref{u}. 


\begin{tiny}
\be
\label{Wh31}  
\begin{aligned}
V_{\mathrm{Wh}(3_1),2}(t,q) = &
  1 + (-2 - 2 q^{-2} + 2 q^{-1} + 4 q - 4 q^2 + 4 q^3 - 2 q^4 - 2 q^7 + 2 q^8 - 
  2 q^{10} + 2 q^{11} + 2 q^{15} - 2 q^{16} + 2 q^{17}  \\ & - 2 q^{18}
 - 2 q^{20} + 
    4 q^{22} - 2 q^{23}) u + (2 + 2 q^{-2} - 2 q^{-1} - 4 q + 2 q^2 - 4 q^3 + 
    2 q^4 + 4 q^5 - 2 q^6 + 4 q^7 \\ & - 4 q^8 + 4 q^9   - 6 q^{10}
+ 2 q^{13} - 
    2 q^{14} + 2 q^{15} + 2 q^{18} - 2 q^{19} + 2 q^{20} - 4 q^{21} + 2 q^{22}) u^2,
\\    
V_{\mathrm{Wh}(\overline{3_1}),2}(t,q) = &
1 + (-2 q^{-26} + 4 q^{-25} - 2 q^{-23} - 2 q^{-21} + 2 q^{-20} - 2 q^{-19}
+ 2 q^{-18} + 2 q^{-14} - 2 q^{-13} + 2 q^{-11} \\ & - 2 q^{-10} + 2 q^{-9}
- 4 q^{-8} + 4 q^{-7} - 6 q^{-6} + 4 q^{-5} - 4 q^{-4} + 4 q^{-3}) u
+ (- 2 q^{-25} + 4 q^{-24} - 2 q^{-23} \\ & + 2 q^{-22} - 2 q^{-21} - 2 q^{-18}
+ 2 q^{-17} - 2 q^{-16} + 6 q^{-13} - 4 q^{-12} + 4 q^{-11} - 4 q^{-10}
+ 2 q^{-9} - 2 q^{-8} \\ & - 2 q^{-7} + 4 q^{-6} - 4 q^{-5} 
+ 2 q^{-4} - 2 q^{-3} + 2 q^{-2}) u^2
\end{aligned}
\ee
\end{tiny}
and for fun,

\begin{tiny}
\be
\label{Wh41}
\begin{aligned}
V_{\mathrm{Wh}(4_1),2}(t,q) = &
1 + (-14 - 2 q^{-18} + 4 q^{-17} + 2 q^{-16} - 6 q^{-15} - 4 q^{-14}
+ 6 q^{-13}
+ 6 q^{-12} - 4 q^{-11} - 8 q^{-10} + 4 q^{-9} \\ & + 8 q^{-8} - 4 q^{-7}
- 10 q^{-6} + 16 q^{-5} - 16 q^{-4} + 20 q^{-3} - 22 q^{-2} + 18 q^{-1}
+ 8 q + 2 q^2 - 4 q^3 - 6 q^4 + 8 q^5 \\ & + 4 q^6 - 8 q^7 - 4 q^8 + 6 q^9
+ 6 q^{10} - 4 q^{11} - 6 q^{12} + 2 q^{13} + 4 q^{14} - 2 q^{15}) u
+ (24 - 2 q^{-17} + 4 q^{-16} \\ & - 2 q^{-14} - 4 q^{-13} + 4 q^{-12}
+ 2 q^{-11} + 4 q^{-10} - 16 q^{-9} + 10 q^{-8} + 8 q^{-7} - 10 q^{-6}
- 4 q^{-5} + 20 q^{-4} \\ & - 28 q^{-3} + 28 q^{-2} - 28 q^{-1} - 16 q + 
    2 q^2 + 12 q^3 - 8 q^4 - 10 q^5 + 16 q^6 - 4 q^7 - 2 q^8 - 
    4 q^9 + 4 q^{10} + 2 q^{11} \\ & - 4 q^{13} + 2 q^{14}) u^2
\end{aligned}
\ee
\end{tiny}
with $u$ as in~\eqref{u}.


\section{Values for  (2,1)-parallel of knots}

We now give values of the $V_{K,1}$, $V_{K(2,1),1}$ and $V_{K,2}$ for some sample
knots to explicitly confirm Equation~\eqref{relV1V2}.


\begin{tiny}
\be
\label{v31}
\begin{aligned}
V_{3_1,1}(t,q) = &
1 + (q + q^3) u + q^2 u^2,
\\  
V_{3_1,2}(t,q) = &
1 + (q + 2 q^3 - q^4 + q^5 - q^6) u + (q^2 + q^4 - q^5) u^2,
\\  
V_{3_1(2,1),1}(t,q) = &
 1 + (q^{-3} + 2 q + 3 q^3 + 2 q^5 + q^7 - q^{13}) u + (3 + 3 q^{-2} + 
 6 q^2 + 4 q^4 + 4 q^6 + 2 q^8 - 2 q^{10} - 2 q^{12}) u^2 \\ &
 + (3 q^{-1} + 3 q + q^3 + q^5 + q^7 - q^{11}) u^3 + u^4
\end{aligned}  
\ee
\end{tiny}
   
\begin{tiny}
\be
\label{v41}
\begin{aligned}
V_{4_1,1}(t,q) = &
1 + (- q^{-1} - q) u + u^2,
\\  
V_{4_1,2}(t,q) = &
1 + (2 - q^{-3} + q^{-2} - 2 q^{-1} - 2 q + q^2 - q^3) u
+ (1 + q^{-2} - q^{-1} - q + q^2) u^2, 
\\  
V_{4_1(2,1),1}(t,q) = &
1 + (- q^{-7} + q^{-5} - 3 q^{-3} - 3 q^{-1} - q^3 - q^7) u
+ (2 + 2 q^{-6} + 5 q^{-4} + q^{-2} + 2 q^4 + 2 q^6) u^2 \\ &
+ (q^{-5} + 3 q^{-3} + 3 q^{-1} + q^5) u^3 + q^{-2} u^4  
\end{aligned}  
\ee
\end{tiny}
   
\begin{tiny}
\be
\label{v61}
\begin{aligned}
V_{6_1,1}(t,q) = &
1 + (- q^{-1} - 2 q - q^3) u + (3 + q^2) u^2,
\\  
V_{6_1,2}(t,q) = &
 1 + (1 - q^{-3} - 2 q^{-1} - 2 q + 2 q^2 - 2 q^3 + 2 q^4 - 2 q^5 + q^6 - 
 q^7) u \\ & 
 + (5 + 3 q^{-2} - 2 q^{-1} - 5 q + 3 q^2 - 3 q^3 + 3 q^4 - q^5 + q^6) u^2, 
\\  
V_{6_1(2,1),1}(t,q) = &
1 + (- q^{-7} + 2 q^{-5} - 5 q^{-3} - 8 q^{-1} - 2 q - q^{11} - q^{15}) u
+ (15 + 6 q^{-6} + 17 q^{-4} + 16 q^{-2} - 2 q^2 - 4 q^4 + 4 q^8 \\ & + 4 q^{10} + 
2 q^{12} + 2 q^{14}) u^2
+ (3 q^{-5} + 10 q^{-3} + 15 q^{-1} + 3 q - 2 q^3 + 2 q^9 + q^{13}) u^3
+ (1 + 3 q^{-2}) u^4  
\end{aligned}  
\ee
\end{tiny}

\begin{tiny}
\be
\label{v62}
\begin{aligned}
V_{6_2,1}(t,q) = &
1 + (- q^{-1} - q^5) u + (1 - q^2 - q^4) u^2 + (q + q^3) u^3 + q^2 u^4,
\\  
V_{6_2,2}(t,q) = &
1 + (1 - q^{-3} + q^{-2} - 2 q^{-1} - q + q^2 - q^5 + q^6 - 2 q^7 + 2 q^8 - q^9) u
+ (-1 + q^{-2} - q^{-1} + 2 q - 4 q^2 + 5 q^3 \\ &
- 6 q^4 + 5 q^5 - 3 q^6 + 2 q^7 - q^8) u^2
+ (-1 + q^{-1} + 3 q - 3 q^2 + 4 q^3 - 4 q^4 + 3 q^5 - 2 q^6 + q^7) u^3 \\ &
+ (1 - q + 2 q^2 - 2 q^3 + 2 q^4 - 2 q^5 + q^6) u^4, 
\\  
V_{6_2(2,1),1}(t,q) = &
 1 + (q^{-9} - 3 q^{-7} + q^{-3} - 6 q^{-1} + 4 q^3 - 2 q^5 - 2 q^7 - q^{13} + 
 q^{17} - q^{19}) u
 + (7 q^{-8} - q^{-6} - 5 q^{-4} - q^{-2} - 8 q^2 \\ & 
 - 12 q^4 - q^6 + q^8 + 2 q^{10} + 2 q^{12} - 2 q^{18}) u^2
 + (21 q^{-7} + 20 q^{-5} + 6 q^{-3} + 10 q^{-1} + 18 q + 14 q^3 + 18 q^5
 + 9 q^7 \\ & - 5 q^9 + 2 q^{11} - q^{13} + 5 q^{15} + 3 q^{17}) u^3
 + (64 + 35 q^{-6} + 62 q^{-4} + 64 q^{-2} + 63 q^2 + 41 q^4 + 22 q^6
 - 12 q^8 \\ & - 16 q^{10} + 16 q^{14} + 12 q^{16}) u^4
 + (35 q^{-5} + 73 q^{-3} + 74 q^{-1} + 51 q + 29 q^3 + 18 q^5 - 13 q^9
 - 11 q^{11} + 11 q^{13} \\ & + 13 q^{15}) u^5
 + (25 + 21 q^{-4} + 39 q^{-2} + 6 q^2 + 7 q^4 - 6 q^{10} + 6 q^{14}) u^6
 + (7 q^{-3} + 8 q^{-1} + q^3 - q^{11} + q^{13}) u^7 + q^{-2} u^8  
\end{aligned}  
\ee
\end{tiny}
and finally, 
\begin{tiny}
\be
\label{v84}
\begin{aligned}
V_{8_4,1}(t,q) = &
1 + (- q^{-3} - 2 q^{-1} - q - q^3 - q^5) u
+ (3 + q^{-2} + q^4) u^2
+ (q^{-1} + 6 q + 5 q^3) u^3 + (1 + 3 q^2) u^4,
\\  
V_{8_4,2}(t,q) = &  
1 + (2 - q^{-7} + q^{-6} - 2 q^{-5} + q^{-4} - 2 q^{-3} + q^{-2}
- 3 q^{-1} - 2 q + 3 q^2 - 2 q^3 + q^4 - 2 q^5 - q^7 + q^8 - q^9) u \\ &
+ (8 + q^{-6} - q^{-5} + 3 q^{-4} - 2 q^{-3} + 3 q^{-2} - 5 q^{-1}
- 7 q + 8 q^2 - 10 q^3 + 8 q^4 - 5 q^5 + 6 q^6 - 3 q^7 + q^8) u^2 \\ & 
+ (-20 + q^{-5} - q^{-4} + 7 q^{-3} - 7 q^{-2} + 16 q^{-1} + 24 q
- 27 q^2 + 27 q^3 - 22 q^4 + 18 q^5 - 9 q^6 + 5 q^7) u^3 \\ & 
+ (6 + q^{-4} - q^{-3} + 4 q^{-2} - 4 q^{-1} - 8 q + 8 q^2 - 7 q^3 + 7 q^4
- 5 q^5 + 3 q^6) u^4, 
\\
V_{8_4(2,1),1}(t,q) = &
1 + (-q^{-15} - 5 q^{-7} - 3 q^{-3} - 11 q^{-1} + 2 q + 6 q^3 - 6 q^5 - 
4 q^7 + 2 q^9 - 2 q^{11} - q^{13} - q^{19}) u
+ (14 + 2 q^{-14} \\ & + 2 q^{-12} + 11 q^{-10} + 24 q^{-8} - 2 q^{-6}
- 9 q^{-4} + 14 q^{-2} - 6 q^2 - 11 q^4 + 17 q^8 + 8 q^{10} + 8 q^{12}
+ 2 q^{14} - 2 q^{16} \\ & + 2 q^{18}) u^2
+ (5 q^{-13} + 8 q^{-11} + 51 q^{-9} + 127 q^{-7} + 124 q^{-5} + 89 q^{-3}
+ 63 q^{-1} + 56 q + 46 q^3 + 49 q^5 + 52 q^7 \\ & + 27 q^9 + 41 q^{11}
+ 27 q^{13} + 22 q^{15} + 21 q^{17}) u^3
+ (141 + 12 q^{-12} + 28 q^{-10} + 107 q^{-8} + 259 q^{-6} + 330 q^{-4} \\ & 
+ 255 q^{-2} + 84 q^2 + 47 q^4 + 90 q^6 + 44 q^8 + 44 q^{10} + 60 q^{12}
+ 60 q^{14} + 44 q^{16}) u^4
+ (13 q^{-11} + 24 q^{-9} + 89 q^{-7} \\ & + 237 q^{-5} + 309 q^{-3} + 215 q^{-1}
+ 70 q + 5 q^3 + 39 q^5 + 41 q^7 + 20 q^9 + 21 q^{11} + 44 q^{13} + 41 q^{15}) u^5
+ (59 \\ & + 6 q^{-10} + 6 q^{-8} + 39 q^{-6} + 114 q^{-4} + 132 q^{-2}
- 11 q^2 + 9 q^4 + 6 q^6 + 12 q^8 + 6 q^{12} + 18 q^{14}) u^6 \\ &
+ (q^{-9} + 10 q^{-5} + 28 q^{-3} + 23 q^{-1} - 2 q + q^5 + 2 q^9 - 2 q^{11}
+ 3 q^{13}) u^7 + (q^{-4} + 3 q^{-2}) u^8  
\end{aligned}  
\ee
\end{tiny}

Keep in mind that the genus of the $(2,1)$-parallel of $K$ is twice the genus of
$K$, and that the knots $3_1$, $4_1$, $6_1$, $6_2$ and $8_4$ have genus
$1, 1, 1, 2, 2$, hence we expect (and we find) that the $V_1$-polynomial of their
$(2,1)$-parallel to have $u$-degree $2, 2, 2, 4, 4$ confirming the equality
in~\eqref{degVn} for the $(2,1)$-parallels of $3_1$, $4_1$, $6_1$, $6_2$ and $8_4$.



\bibliographystyle{hamsalpha}
\bibliography{biblio}

\newcommand{\etalchar}[1]{$^{#1}$}
\providecommand{\bysame}{\leavevmode\hbox to3em{\hrulefill}\thinspace}
\providecommand{\href}[2]{#2}
\providecommand{\eprint}{\begingroup \urlstyle{rm}\Url}
\begin{thebibliography}{BNvdV24}

\bibitem[AS02]{AS:pointed}
Nicol\'{a}s Andruskiewitsch and Hans-J\"{u}rgen Schneider, \emph{Pointed {H}opf algebras}, New directions in {H}opf algebras, Math. Sci. Res. Inst. Publ., vol.~43, Cambridge Univ. Press, Cambridge, 2002, pp.~1--68.

\bibitem[BC17]{JC:tensornet}
Jacob Bridgeman and Christopher Chubb, \emph{Hand-waving and interpretive dance: an introductory course on tensor networks}, J. Phys. A \textbf{50} (2017), no.~22, 223001, 61.

\bibitem[BNvdV]{BNV:PG}
Dror Bar-Natan and Roland van~der Veen, \emph{Perturbed {G}aussian generating functions for universal knot invariants}, Preprint 2021, \href{https://arxiv.org/abs/2109.02057}{arXiv:2109.02057}.

\bibitem[BNvdV24]{BNV:API}
\bysame, \emph{A perturbed-{A}lexander invariant}, Quantum Topol. \textbf{15} (2024), no.~3, 449--472.

\bibitem[Bur20]{Burton:knots}
Benjamin Burton, \emph{The next 350 million knots}, 36th {I}nternational {S}ymposium on {C}omputational {G}eometry, LIPIcs. Leibniz Int. Proc. Inform., vol. 164, Schloss Dagstuhl. Leibniz-Zent. Inform., Wadern, 2020, pp.~Art. No. 25, 17.

\bibitem[CDGW]{snappy}
Marc Culler, Nathan Dunfield, Matthias Goerner, and Jeffrey Weeks, \emph{Snap{P}y, a computer program for studying the geometry and topology of {$3$}-manifolds}, Available at \url{http://snappy.computop.org}.

\bibitem[DFJ12]{Dunfield:twisted}
Nathan Dunfield, Stefan Friedl, and Nicholas Jackson, \emph{Twisted {A}lexander polynomials of hyperbolic knots}, Exp. Math. \textbf{21} (2012), no.~4, 329--352.

\bibitem[DGST10]{Dunfield:mutation}
Nathan Dunfield, Stavros Garoufalidis, Alexander Shumakovitch, and Morwen Thistlethwaite, \emph{Behavior of knot invariants under genus 2 mutation}, New York J. Math. \textbf{16} (2010), 99--123.

\bibitem[Gab86]{Gabai}
David Gabai, \emph{Genera of the arborescent links}, Mem. Amer. Math. Soc. \textbf{59} (1986), no.~339, i--viii and 1--98.

\bibitem[Gar04]{Ga:whitehead}
Stavros Garoufalidis, \emph{Whitehead doubling persists}, Algebr. Geom. Topol. \textbf{4} (2004), 935--942.

\bibitem[GHK{\etalchar{+}}a]{GHKST:link}
Stavros Garoufaldis, Matthew Harper, Ben-Michael Kohli, Jiebo Song, and Guillaume Tahar, \emph{Extending knot polynomials of braided {H}opf algebras to links}, Preprint 2025, \href{https://arxiv.org/abs/2505.01398}{arXiv:2505.01398}.

\bibitem[GHK{\etalchar{+}}b]{Vn}
Stavros Garoufalidis, Matthew Harper, Rinat Kashaev, Ben-Michael Kohli, and Emmanuel Wagner, \emph{On the colored {L}inks--{G}ould polynomial}, Preprint 2025, \href{https://arxiv.org/abs/2509.10911}{arXiv:2509.10911}.

\bibitem[GK26]{GK:multi}
Stavros Garoufalidis and Rinat Kashaev, \emph{Multivariable {K}not {P}olynomials from {B}raided {H}opf {A}lgebras with {A}utomorphisms}, Publ. Res. Inst. Math. Sci. \textbf{62} (2026), no.~1, 75--114.

\bibitem[GKK{\etalchar{+}}]{LG-V1}
Stavros Garoufalidis, Rinat Kashaev, Ben-Michael Kohli, Jiebo Song, and Guillaume Tahar, \emph{Skein theory for the {L}inks-{G}ould polynomial}, to appear in the Journal of the London Math.Soc.

\bibitem[GL24]{GS:VnData}
Stavros Garoufalidis and Shana Li, \emph{Values of {$V_n$}-polynomials of knots}, \url{https://dataverse.harvard.edu/dataset.xhtml?persistentId=doi:10.7910/DVN/XE4TOF}, 2024.

\bibitem[HTW98]{HTW}
Jim Hoste, Morwen Thistlethwaite, and Jeff Weeks, \emph{The first 1,701,936 knots}, Math. Intelligencer \textbf{20} (1998), no.~4, 33--48.

\bibitem[HW18]{Hedden:geography}
Matthew Hedden and Liam Watson, \emph{On the geography and botany of knot {F}loer homology}, Selecta Math. (N.S.) \textbf{24} (2018), no.~2, 997--1037.

\bibitem[Jab]{Jablan}
Slavik Jablan, \emph{Tables of quasi-alternating knots with at most 12 crossings}, Preprint 2014, \href{https://arxiv.org/abs/1404.4965}{ arXiv:1404.4965}.

\bibitem[Jon87]{Jones}
Vaughan Jones, \emph{Hecke algebra representations of braid groups and link polynomials}, Ann. of Math. (2) \textbf{126} (1987), no.~2, 335--388.

\bibitem[Kas21]{Ka:longknots}
Rinat Kashaev, \emph{Invariants of long knots}, Representation theory, mathematical physics, and integrable systems, Progr. Math., vol. 340, Birkh\"{a}user/Springer, Cham, [2021] \copyright 2021, pp.~431--451.

\bibitem[Kno]{knotatlas}
KnotAtlas, Available at \url{http://katlas.org}.

\bibitem[KT]{Kohli-Tahar}
Ben-Michael Kohli and Guillaume Tahar, \emph{A lower bound for the genus of a knot using the {L}inks--{G}ould invariant}, Preprint 2023, \href{https://arxiv.org/abs/2310.15617}{arXiv:2310.15617}.

\bibitem[KT57]{KT}
Shin'ichi Kinoshita and Hidetaka Terasaka, \emph{On unions of knots}, Osaka Math. J. \textbf{9} (1957), 131--153.

\bibitem[KWZ]{KWZ}
Artem Kotelskiy, Liam Watson, and Claudius Zibrowius, \emph{Immersed curves in {K}hovanov homology}, Preprint 2019, \href{https://arxiv.org/abs/1910.14584}{ arXiv:1910.14584}.

\bibitem[LG92]{Links-Gould}
Jon Links and Mark Gould, \emph{Two variable link polynomials from quantum supergroups}, Lett. Math. Phys. \textbf{26} (1992), no.~3, 187--198.

\bibitem[Li]{Shana:FPT-RT}
Shana~Yunsheng Li, \emph{Fixed-parameter tractable computation of knot polynomials of reshetikhin--turaev type via tensor networks}, In preparation.

\bibitem[LSW97]{LSR:NPtensor}
Chi-Chung Lam, P.~Sadayappan, and Rephael Wenger, \emph{On optimizing a class of multi-dimensional loops with reduction for parallel execution}, Parallel Process. Lett. \textbf{7} (1997), no.~2, 157--168.

\bibitem[MO08]{Manolescu:thin}
Ciprian Manolescu and Peter Ozsv\'{a}th, \emph{On the {K}hovanov and knot {F}loer homologies of quasi-alternating links}, Proceedings of {G}\"{o}kova {G}eometry-{T}opology {C}onference 2007, G\"{o}kova Geometry/Topology Conference (GGT), G\"{o}kova, 2008, pp.~60--81.

\bibitem[Mor86]{Morton:cg}
Hugh Morton, \emph{Seifert circles and knot polynomials}, Math. Proc. Cambridge Philos. Soc. \textbf{99} (1986), no.~1, 107--109.

\bibitem[MST]{MST:QES}
Scott Morrison, Noah Snyder, and Dylan Thurston, \emph{Towards the quantum exceptional series}, Preprint 2024, \href{https://arxiv.org/abs/2402.03637}{arXiv:2402.03637}.

\bibitem[Mur58]{Murasugi}
Kunio Murasugi, \emph{On the genus of the alternating knot. {I}, {II}}, J. Math. Soc. Japan \textbf{10} (1958), 94--105, 235--248.

\bibitem[ORS13]{ORS:ODD}
Peter Ozsv\'{a}th, Jacob Rasmussen, and Zolt\'{a}n Szab\'{o}, \emph{Odd {K}hovanov homology}, Algebr. Geom. Topol. \textbf{13} (2013), no.~3, 1465--1488.

\bibitem[OS04]{OS:genusbounds}
Peter Ozsv\'{a}th and Zolt\'{a}n Szab\'{o}, \emph{Knot {F}loer homology, genus bounds, and mutation}, Topology Appl. \textbf{141} (2004), no.~1-3, 59--85.

\bibitem[OS05]{OS:quasi}
\bysame, \emph{On the {H}eegaard {F}loer homology of branched double-covers}, Adv. Math. \textbf{194} (2005), no.~1, 1--33.

\bibitem[OS06]{OS}
\bysame, \emph{An introduction to {H}eegaard {F}loer homology}, Floer homology, gauge theory, and low-dimensional topology, Clay Math. Proc., vol.~5, Amer. Math. Soc., Providence, RI, 2006, pp.~3--27.

\bibitem[RT90]{RT:ribbon}
Nikolai Reshetikhin and Vladimir Turaev, \emph{Ribbon graphs and their invariants derived from quantum groups}, Comm. Math. Phys. \textbf{127} (1990), no.~1, 1--26.

\bibitem[SS22]{SS:thick}
Andr\'as~I. Stipsicz and Zolt\'an Szab\'o, \emph{A note on thickness of knots}, Gauge theory and low-dimensional topology---progress and interaction, Open Book Ser., vol.~5, Math. Sci. Publ., Berkeley, CA, 2022, pp.~299--308.

\bibitem[Sto]{Stoimenow}
Alexander Stoimenow, \emph{Tables of mutant knots}, Available at \url{http://stoimenov.net/stoimeno/homepage/ptab/index.html}.

\bibitem[Tur94]{Tu:book}
Vladimir Turaev, \emph{Quantum invariants of knots and 3-manifolds}, de Gruyter Studies in Mathematics, vol.~18, Walter de Gruyter \& Co., Berlin, 1994.

\end{thebibliography}
\end{document}